\documentclass[12pt,psamsfonts,leqno,oneside,letterpaper]{amsart}
\usepackage[dvips,text={6.5truein,9truein},left=1truein,top=1truein]{geometry}
\usepackage{amssymb,amsmath,amscd,enumerate,
}
\usepackage[pdftex]{graphicx}
\usepackage{url}

\usepackage[colorlinks,linkcolor=blue,citecolor=blue,pdfstartview=FitH]{hyperref}
\input xy
\xyoption{all}
\SelectTips{cm}{12}

\usepackage{color}
\usepackage{mathrsfs}
\newcommand\invisiblecomment[1]{\empty}

\newcommand\missingref[1]{\empty}
\newcommand\tinymissingref[1]{\empty}
\newcommand\abstractcomment[1]{\empty}

\theoremstyle{definition}
\newtheorem{para}{}[section]

\newtheorem{remark}[para]{Remark}
\newtheorem{reformulation}[para]{Reformulation}
\newtheorem{remarks}[para]{Remarks}
\newtheorem{notation}[para]{Notation}

\newtheorem{convention}[para]{Convention}
\newtheorem{definition}[para]{Definition}
\newtheorem{definitions}[para]{Definitions}
\newtheorem{definitionnotation}[para]{Definition and Notation}
\newtheorem{remarksnotation}[para]{Remarks and Notation}

\newtheorem{remarknotation}[para]{Remark and Notation}
\newtheorem{notationremark}[para]{Notation and Remark}
\newtheorem{notationreviewremarks}[para]{Notation, Review and Remarks}
\newtheorem{definitionremark}[para]{Definition and Remark}
\newtheorem{definitionsremarks}[para]{Definitions and Remarks}
\newtheorem{notationremarks}[para]{Notation and Remarks}
\newtheorem{definitionsnotation}[para]{Definitions and Notation}

\newtheorem{reviewdefinition}[para]{Review and Definition}

\newtheorem{remarksdefinition}[para]{Remarks and Definition}
\newtheorem{definitionnotationremarks}[para]{Definition, Notation and Remarks}
\newtheorem{definitionsnotationremarks}[para]{Definitions, Notation and Remarks}

\newcommand\Alternatives{\begin{enumerate}[(i)]}
\newcommand\EndAlternatives{\end{enumerate}}
\newcommand\Conditions{\begin{enumerate}[(1)]}
\newcommand\EndConditions{\end{enumerate}}

\theoremstyle{plain}
\newtheorem{theorem}[para]{Theorem}
\newtheorem{lemma}[para]{Lemma}
\newtheorem{proposition}[para]{Proposition}

\newtheorem{corollary}[para]{Corollary}
\newtheorem{conjecture}[para]{Conjecture}
\newtheorem*{introtheorem}{Theorem}

\newtheorem{claim}[equation]{}
\numberwithin{equation}{para}
\numberwithin{figure}{section}
\numberwithin{specialremark}{para}

\newcommand\Number{\begin{para}}
\newcommand\EndNumber{\end{para}}
\newcommand\Definition{\begin{definition}}
\newcommand\EndDefinition{\end{definition}}
\newcommand\Definitions{\begin{definitions}}
\newcommand\DefinitionsNotation{\begin{definitionsnotation}}
\newcommand\DefinitionNotation{\begin{definitionnotation}}
\newcommand\RemarksNotation{\begin{remarksnotation}}
\newcommand\Reformulation{\begin{reformulation}}
\newcommand\EndRemarksNotation{\end{remarksnotation}}
\newcommand\EndReformulation{\end{reformulation}}
\newcommand\RemarkNotation{\begin{remarknotation}}
\newcommand\EndRemarkNotation{\end{remarknotation}}
\newcommand\NotationRemark{\begin{notationremark}}
\newcommand\EndDefinitionNotationRemarks{\end{definitionnotationremarks}}
\newcommand\NotationReviewRemarks{\begin{notationreviewremarks}}
\newcommand\DefinitionRemark{\begin{definitionremark}}
\newcommand\RemarksDefinition{\begin{remarksdefinition}}
\newcommand\DefinitionsRemarks{\begin{definitionsremarks}}
\newcommand\DefinitionNotationRemarks{\begin{definitionnotationremarks}}
\newcommand\DefinitionsNotationRemarks{\begin{definitionsnotationremarks}}
\newcommand\EndDefinitionsNotationRemarks{\end{definitionsnotationremarks}}

\newcommand\NotationRemarks{\begin{notationremarks}}
\newcommand\EndNotationRemark{\end{notationremark}}
\newcommand\EndNotationReviewRemarks{\end{notationreviewremarks}}
\newcommand\EndDefinitionRemark{\end{definitionremark}}
\newcommand\EndRemarksDefinition{\end{remarksdefinition}}
\newcommand\EndDefinitionsRemarks{\end{definitionsremarks}}
\newcommand\EndNotationRemarks{\end{notationremarks}}
\newcommand\EndDefinitionsNotation{\end{definitionsnotation}}
\newcommand\EndDefinitionNotation{\end{definitionnotation}}
\newcommand\ReviewDefinition{\begin{reviewdefinition}}
\newcommand\EndReviewDefinition{\end{reviewdefinition}}
\newcommand\EndDefinitions{\end{definitions}}
\newcommand\Theorem{\begin{theorem}}
\newcommand\EndTheorem{\end{theorem}}
\newcommand\Conjecture{\begin{conjecture}}
\newcommand\EndConjecture{\end{conjecture}}
\newcommand\Remark{\begin{remark}}
\newcommand\EndRemark{\end{remark}}
\newcommand\Remarks{\begin{remarks}}
\newcommand\EndRemarks{\end{remarks}}
\newcommand\Convention{\begin{convention}}
\newcommand\EndConvention{\end{convention}}
\newcommand\Notation{\begin{notation}}
\newcommand\EndNotation{\end{notation}}
\newcommand\Lemma{\begin{lemma}}
\newcommand\EndLemma{\end{lemma}}
\newcommand\Proposition{\begin{proposition}}
\newcommand\EndProposition{\end{proposition}}
\newcommand\Corollary{\begin{corollary}}
\newcommand\EndCorollary{\end{corollary}}
\newcommand\Claim{\begin{claim}}
\newcommand\EndClaim{\end{claim}}
\newcommand\Proof{\begin{proof}}
\newcommand\EndProof{\end{proof}}
\newcommand\Equation{\begin{equation}}
\newcommand\EndEquation{\end{equation}}

\newcommand\Bullets{\begin{itemize}}
\newcommand\EndBullets{\end{itemize}}

\newcommand\dotsystem{dot system}

\newcommand\realcalH{\mathcal{H}}
\newcommand\calH{\mathscr{H}}
\newcommand\redPi{\Pi}
\newcommand\redL{L}
\newcommand\redXi{\Xi}
\newcommand\magentapacking{discrete set}
\newcommand\redE{E}
\newcommand\rede{e}
\newcommand\redcalH{\calH}
\newcommand\redrealcalH{\realcalH}

\renewcommand\epsilon{\varepsilon}

\newcommand\RR{{\bf R}}

\newcommand\hatgamma{\widehat\gamma}

\newcommand\XPS{\calX_{P,\scrS}}
\newcommand\XPtS{\calX_{P,\tS}}
\newcommand\maybeL{\mathcal V}

\newcommand\maybescrK{\scrK}

\newcommand\Vlens{V_{\rm lens}}
\newcommand\Vcone{V_{\rm cone}}
\newcommand\Wlens{W_{\rm lens}}
\newcommand\Wcone{W_{\rm cone}}

\newcommand\inter{\mathop{\rm int}}

\newcommand\Mthick{M_{\rm thick}}
\newcommand\Mthin{M_{\rm thin}}

\newcommand\tH{\widetilde H}

\newcommand\wasx{r_1}
\newcommand\wasy{r_2}
\newcommand\wasz{D}

\newcommand\maybecalD{\calD}
\newcommand\willbep{p}
\newcommand\redu{u}
\newcommand\redU{U}
\newcommand\redtau{\tau}
\newcommand\tredtau{T}
\newcommand\tryepsilon{\epsilon}
\newcommand\willbeX{X}
\newcommand\willbeS{S}
\newcommand\willbeP{P}
\newcommand\willbei{i}
\newcommand\willbefinitevolume{finite-volume}
\newcommand\willbesemifree{semifree}
\newcommand\hatp{{\widehat p}}
\newcommand\willbehatp{{\widehat p}}
\newcommand\willbepolytope{polyhedron}
\newcommand\willbelambdanought{\lambda_0}
\newcommand\willbepolytopal{polyhedral}
\newcommand\willbepolyhedral{polyhedral}

\newcommand\willbepolytopein{polyhedron in\ }
\newcommand\willbepolyhedra{polyhedra}

\newcommand\maybeU{T}
\newcommand\maybeepsilon{\epsilon}
\newcommand\maybeb{b}

\newcommand\tS{{\widetilde S}}

\newcommand\calc{{\mathcal C}}
\newcommand\calD{{\mathcal D}}

\newcommand\cala{{\mathcal A}}

\newcommand\calf{{\mathcal F}}

\newcommand\realcalX{{\mathcal X}}
\newcommand\calX{X}

\newcommand\anoughtS{{\cala^0_S}}
\newcommand\aoneS{{\cala^1_S}}
\newcommand\atwoS{{\cala^2_S}}
\newcommand\athreeS{{\cala^3_S}}

\newcommand\arccosh{{\rm arccosh}}

\newcommand\arcsech{{\rm arcsech}}

\newcommand\sech{{\rm sech}}

\newcommand\arcsec{{\rm arcsec}}
\newcommand\ZZ{{\bf Z}}

\newcommand\HH{{\bf H}}
\newcommand\EE{{\bf E}}

\newcommand\cals{{\mathscr S}}

\newcommand\scrD{{\mathscr D}}
\newcommand\scrG{{\mathcal G}}

\newcommand\scrK{{\mathscr K}}

\newcommand\scrP{{\mathscr P}}
\newcommand\redscrPST{\scrP_{S,\redcalT}}

\newcommand\scrS{{\mathscr S}}
\newcommand\scrT{{\mathscr T}}
\newcommand\scrV{{\mathscr V}}

\newcommand\scrv{{\mathscr V}}

\newcommand\dist{\mathop{\rm dist}}

\newcommand\nbhd{\mathop{\rm nbhd}}
\newcommand\clnbhd{\nbhd}

\newcommand\density{{\rm density}}

\newcommand\vol{\mathop{\rm vol}}
\newcommand\area{\mathop{\rm area}}

\newcommand\voct{V_{\rm oct}}

\newcommand\length{\mathop{{\rm length}}}

\newcommand\gpc{generalized polyhedral complex}

\newcommand\rank{\mathop{{\rm rank}}}

\newcommand\image{\mathop{{\rm Im}}}

\newcommand\isomplus{\mathop{{\rm Isom}_+}}
\newcommand\isom{\mathop{{\rm Isom}}}

\newcommand\FF{{\bf F}}

\newcommand\redcalM{\mathcal M}
\newcommand\redcalT{\mathcal T}

\newcommand\mayber{r}

\DeclareFontFamily{U}{rcjhbltx}{}
\DeclareFontShape{U}{rcjhbltx}{m}{n}{<->rcjhbltx}{}
\DeclareSymbolFont{hebrewletters}{U}{rcjhbltx}{m}{n}

\let\aleph\relax\let\beth\relax
\let\gimel\relax\let\daleth\relax

\DeclareMathSymbol{\aleph}{\mathord}{hebrewletters}{39}
\DeclareMathSymbol{\beth}{\mathord}{hebrewletters}{98}
\DeclareMathSymbol{\gimel}{\mathord}{hebrewletters}{103}
\DeclareMathSymbol{\daleth}{\mathord}{hebrewletters}{100}

\DeclareMathSymbol{\lamed}{\mathord}{hebrewletters}{108}
\DeclareMathSymbol{\mem}{\mathord}{hebrewletters}{109}
\DeclareMathSymbol{\ayin}{\mathord}{hebrewletters}{96}
\DeclareMathSymbol{\tsadi}{\mathord}{hebrewletters}{118}
\DeclareMathSymbol{\qof}{\mathord}{hebrewletters}{114}
\DeclareMathSymbol{\shin}{\mathord}{hebrewletters}{152}

\usepackage[utf8]{inputenc}

\begin{document}

\title{The ratio of homology rank to hyperbolic volume, I}
\author{Rosemary K. Guzman}
\address{Department of Mathematics
\\
University of Illinois\\
1409 W. Green St.\\
Urbana, IL 61801}
\email{rguzma1@illinois.edu}

\author{Peter B. Shalen}
\address{Department of Mathematics, Statistics, and Computer Science
(M/C 249)\\
University of Illinois at Chicago\\
851 S. Morgan St.\\
Chicago, IL 60607-7045}
\email{petershalen@gmail.com}


\maketitle

\begin{abstract}
We show that for
every finite-volume 
hyperbolic $3$-manifold $M$ and every prime $p$ we
have 
$\dim H_1(M;\FF_p)< 
168.602
\cdot\vol 
(M)
$.
There are slightly stronger estimates if $p=2$ or if $M$ is non-compact.
This improves on a result
proved by Agol, Leininger and Margalit, 
which gave 
the same inequality with a coefficient of
$334.08$
in place of 
$168.602$. It also improves on the analogous result with a
coefficient of about $260$, which could have been obtained 
by combining the
arguments due to Agol, Leininger and Margalit  with a result due to B\"or\"oczky. Our inequality
involving homology rank is deduced from a result about the rank of the
fundamental group: if $M$ is a 
finite-volume 
orientable hyperbolic $3$-manifold such that $\pi_1(M)$ is 
$2$-semifree, 
then
$\rank\pi_1(M)<1+\willbelambdanought\cdot\vol M$, where
$\willbelambdanought$ is a certain constant less than $167.79$.
\end{abstract}

\section{Introduction}

It is a standard consequence of the Margulis Lemma
that for every $n\ge2$ there exists a constant $\lambda>0$ such that for
every finite-volume hyperbolic $n$-manifold $M$ we have
$\rank\pi_1(M)\le\lambda\cdot\vol M$ (where $\vol$ denotes hyperbolic volume).
In particular, there is a constant $\lambda>0$ such that for every
\willbefinitevolume
\ orientable
hyperbolic $3$-manifold and every prime $p$ we have
\Equation\label{in general}
\dim H_1(M;\FF_p)\le\lambda\cdot\vol M.
\EndEquation
According to \cite[Proposition 2.2]{alm}, 
(\ref{in general})
always holds with $\lambda=334.08$. The authors of \cite{alm} have
informed us that they were aware that the result could be improved
using the sphere-packing results of \cite{boroczky}; a straightforward
application of these results shows that (\ref{in general}) holds with
a value of
$\lambda$ 
that is approximately equal to 
$260$.

Theorem \ref{main} of this paper 
includes the assertion
that (\ref{in general})
always holds with a value of $\lambda$ 
just under $168.602$.

Thus we have:
\begin{introtheorem} For any \willbefinitevolume\
orientable hyperbolic $3$-manifold $M$,  and any prime $p$,
we have
$$\dim H_1(M;\FF_p)< 168.602\cdot\vol
( M).
$$
\end{introtheorem}

Theorem \ref{main} also asserts that if one takes $p=2$, or restricts
attention to the case where $M$ is non-compact, then the constant in
the conclusion of the theorem stated above can be improved in the
first place to the right of the decimal point.

Theorem \ref{main} is deduced
from Proposition \ref{pi-one rank}, which gives a bound on the rank
of a
fundamental group under a certain condition. Recall that a group $G$ is
said to be {\it $k$-semifree} for a given positive integer $k$ if each 
subgroup of $G$ having rank at most $k$ is a free product of free
abelian groups.
Proposition
\ref{pi-one rank} asserts that if $M$ is a \willbefinitevolume\
orientable hyperbolic $3$-manifold 
such 
that $\pi_1(M)$ is $2$-\willbesemifree, then
$\rank\pi_1(M)<1+\willbelambdanought\cdot\vol M$,
where $\willbelambdanought$ is a certain constant less than $167.79$.

In the sequels to this paper, by combining the methods of the present paper with a 
variety of deep results
and techniques (including the Four Color Theorem and the
desingularization techniques of
 \cite{acs-singular} and \cite{singular-two}),  we will prove results that significantly improve the estimates
 given by Proposition  \ref{pi-one rank} and Theorem \ref{main}. While
 some of these results involve mild additional
 topological hypotheses, there does appear to be a strict improvement
 of Theorem \ref{main} for the case 
where
$p=2$ and $M$ is closed.

The proof of Proposition \ref{pi-one rank} is a refinement of the
corresponding result
in \cite{alm} (see the final displayed formula in the proof of
\cite[Proposition 2.2]{alm}). The refinement uses a crucial new idea to improve the estimate. In
order to describe this improvement, we shall begin by reviewing the
relevant argument given
in \cite{alm},
using some terms that are explained in the body of the present paper. 

According to \cite[Corollary 4.2]{acs-surgery}, the hypothesis that
the \willbefinitevolume\ orientable $3$-manifold 
$M$ has a $2$-\willbesemifree\ fundamental group implies that $\tryepsilon:=\log3$ is a Margulis
number (see Definitions \ref{marg def} below) for $M$. One fixes a
finite subset $S$ of $M$ which is a maximal
$\tryepsilon$-\magentapacking\ 
contained in the
$\tryepsilon$-thick part of $M$ (see Definitions 
\ref{nets and such one} and \ref{thin def}). 
By a
Voronoi-Dirichlet construction, which we explain in considerable
detail in Section \ref{Voronoi section} below, one associates with the
set $S$ a ``cell complex'' 
 $K_S$ whose underlying space is $M$. (When $M$ is closed, $K_S$ is a
 certain kind of CW complex.) Each (open)
$3$-cell of $K_S$ is the homeomorphic image  under a locally isometric
covering map $q:\HH^3\to M$ of the interior of a
``Voronoi region,''
which is a convex \willbepolytope\  $X$ associated with a point $P$ of
$q^{-1}(S)$, and having $P$ as an interior point. We think of $P$ as
the ``center'' of the Voronoi region $X$.

Now
$\pi_1(M)$ is carried by 
a connected graph $\scrG$ 
which has
$S$ as its vertex set, and has one edge ``dual'' to each $2$-cell of $K_S$
meeting the 
$\epsilon$-thick 
part of $M$. To bound the rank of $\pi_1(M)$ as
a multiple of $\vol M$, it
therefore suffices to bound the first betti number of $\scrG$; this is
done by bounding the number of vertices of $\scrG$ as a
multiple of $\vol M$, and bounding the valence of an arbitrary vertex.
In order to bound the number of vertices of $\scrG$ (or equivalently
the cardinality of $S$) as a
multiple of $\vol M$, it suffices to give a lower bound for the volume
of an arbitrary $3$-cell in $K_S$. In \cite{alm} this is done by
observing that each $3$-cell contains an isometric copy of a ball of
radius $\tryepsilon/2$ in $\HH^3$, and
using the volume of this ball as a lower bound; as mentioned above,
this can be improved using the sphere-packing results of
\cite{boroczky}, and we have incorporated this in our estimates.

Bounding the valence of a vertex of $\scrG$ amounts to bounding, for a
given Voronoi region $X$, the number of
two-dimensional faces of $X$ which meet
the pre-image under $q$ of the $\tryepsilon$-thick part of $M$. In the
present sketch we shall refer to such faces as ``thick faces.''
Any thick face $F$ determines a Voronoi
region $Y$ which ``neighbors'' $X$ along $F$ in the sense that $F=X\cap Y$. 
The neighboring Voronoi region $Y$  
contains the ball $E$ of radius $\tryepsilon/2$ about the
``center'' of $Y$; 
furthermore, the ``centers'' of $X$ and $Y$ are separated by a
distance of at most $2\epsilon$, which implies that 
$E$ is
contained in a ball 
$N_0$
of radius $5\tryepsilon/2$ about the ``center'' $P$ of $X$. This implies an upper bound of
$\lfloor B(5\tryepsilon/2)/B(\tryepsilon/2)-1\rfloor=493$ for the valence of a vertex of $\scrG$, where
$B(r)$ denotes the volume of a ball of radius $r$ in $\HH^3$. This is
the valence bound used in \cite{alm}.

The main new idea in the proof of Proposition \ref{central} involves a somewhat
different approach to bounding the valence of a vertex.
We fix a
constant $R$ between $2\tryepsilon$ and $5\tryepsilon/2$, and denote by $N$ the
closed ball
of radius $R$ about the point $P$. The
strategy is to find a 
relatively 
large constant $c>0$ 
such that, if $F$ is a thick face of $X$, and  $Y$ is the neighboring
Voronoi region determined by $F$,
the volume
of $Y\cap N$ is bounded below by $c$. This will give an upper bound
slightly less than $B(R)/c$ for
the valence of a vertex of $\scrG$.

To carry out this strategy, we choose a
point $\tredtau $ that lies  in the intersectiton of a given thick face $F$ with the preimage under $q$ of the
$\tryepsilon$-thick part of 
$M$. 
The definition of a Voronoi region implies that the 
distance  $D$
from $\tredtau $ to $P$ is equal to the distance from $\tredtau$ to the ``center''
$Q$ of the neighboring Voronoi region
$Y$. Then $Y$ contains the set $Z$
(an ``ice cream cone'') 
which is defined as the convex hull
of the union of $\{\tredtau \}$ with the closed ball of radius $\tryepsilon/2$ about
$Q$, while  $N$  contains 
the closed ball 
of radius $\rho:=R-D$ centered at $\tredtau$. In the notation that is
formalized later in this introduction, this closed ball is denoted by
$\overline{\nbhd_\rho(\tredtau)}$. Since $Y\supset Z$ and
$N\supset\overline{\nbhd_\rho(\tredtau)}$, the set
$Y\cap N$ contains 
$Z\cap \overline{\nbhd_\rho(\tredtau)} $.

One can use elementary hyperbolic geometry to calculate
the volume of 
$Z\cap \overline{\nbhd_\rho(\tredtau)} $ 
in terms of $R$, $\tryepsilon$ and $D$. If we recall that
$\tryepsilon=\log3$, and 
take $R=2\log3+0.15$
(a choice of $R$ that will turn out to minimize our valence bound),
the expression for the
volume of
$Z\cap \overline{\nbhd_\rho(\tredtau)} $ 
becomes a function of the parameter $D$, which can take any
value in the interval $[\tryepsilon/2,\tryepsilon]$. We show that $c:=0.496$ is a lower
bound for this function on the subinterval $[R/2-\tryepsilon/4,\tryepsilon]$. When
$\tryepsilon/2\le D\le R/2-\tryepsilon/4$, it turns out that in the geometric
situation described above, $Y\cap N$ contains the ball of radius
$\tryepsilon/2$ about $Q$; this gives a lower bound of $B(\tryepsilon/2)>c$ for the
volume of $Y\cap N$. Thus $c$ is a lower bound for this volume in all
cases.

This provides an improved upper bound of $314$ for the valence of a
vertex of $\scrG$.
This improved valence bound is the crucial new ingredient in the proof
of Proposition \ref{pi-one rank}.

Theorem \ref{main} is deduced
from Proposition \ref{pi-one rank} by the same method 
as is
used for the
corresponding step in \cite{alm}. 
The constants that appear in Theorem \ref{main}
incorporate
small improvements that are provided by results from
\cite{acs-surgery}, \cite{gmm} and \cite{gs}.

The definition and properties of the  complex $K_S$ are presented in
Section \ref{Voronoi section} of the present paper. This section
contains a good deal of expository material about convex \willbepolyhedra\ and
Voronoi decompositions,
since the foundations of the theories in question do not seem to be  extensively
documented in the literature; our main 
references
are to
\cite{Grunbaum} and
\cite{jason}. In Section \ref{convex section} we develop the quantitative hyperbolic geometry needed to compute
the volume of the set which is denoted 
$Z\cap \overline{\nbhd_\rho(\tredtau)} $ 
in the above
discussion. In Section \ref{central section} we give a general version, for
any Margulis number $\tryepsilon$, of the
argument sketched above for the case $\tryepsilon=\log3$; the result is
embodied in Corollary \ref{quote me}. We regard this corollary, and
Proposition \ref{central} from which it is deduced, as the central
results of this paper, although their statements are more technical
than that of Theorem \ref{main}. In Section \ref{numerical section}, we do the
numerical calculations needed to pass from Corollary \ref{quote me} to
Proposition  \ref{pi-one rank} and Theorem \ref{main}.

We summarize here some conventions that will be used in the body of
the paper.

If $X$ is a group, we will write  $Y\le X$ to mean that $Y$ is a subgroup of
  $X$. 

The symbol ``$\dist$'' will denote the distance in a metric space when it
is clear from the context which metric space is involved. If $p$ is a point of a metric space, and $r$ is a real number, we will denote by $\nbhd_r(p)$ the set of all $x\in X$ such that $\dist(x,p)<r$. Thus $\nbhd_r(p)$  is a neighborhood of $p$ if $r>0$, and is empty if $r\le0$.
 
The isometry
group of a metric space $X$ is denoted $\isom(X)$.
If $X$ is an orientable Riemannian manifold, $\isomplus(X)$ will denote
the orientation-preserving subgroup of  $\isom(X)$.

We thank Ian Agol, Chris Leininger, and Dan Margalit for useful
discussions of their paper \cite{alm}.

\section{The Voronoi complex}\label{Voronoi section}

In this section, after a brief review of facts about convex polyhedra
in $\HH^n$, we introduce the notion of a polyhedral complex in
$\HH^n$, and show how a discrete set in $\HH^n$ defines a polyhedral
complex, its Voronoi complex. We then show how, if $S$ is a finite subset of a
hyperbolic $n$-manifold $M$, the Voronoi complex defined by the pre-image of
$S$ under a locally isometric covering map $\HH^n\to M$ induces
a ``generalized polyhedral complex'' with $M$ as its underlying space,
and an associated dual graph.

\DefinitionsNotationRemarks\label{polyhedra}

Let 
$n$ be a positive integer. A {\it closed convex set} in $\EE^n$
or $\HH^n$ may be defined to be a
non-empty set which is the 
intersection of a collection of closed
half-spaces, or to be a 
non-empty
closed set which contains the line segment joining
any two of its points; 
in the hyperbolic case, the equivalence of the
two definitions follows from \cite[Lemma 1.5]{jason}. 
If $X$ is any closed convex set in $\EE^n$
or $\HH^n$,
 and if $W$ denotes the smallest totally geodesic subspace 
of $\EE^n$ or $\HH^n$ containing $X$, we define 
the {\it dimension} of $X$ to be the dimension of $W$, and we define
the {\it interior} of $X$, denoted $\inter X$, to be its topological
interior  relative to $W$. The {\it boundary} of $X$ is
$\partial X:=X-\inter X$.

A {\it support hyperplane} for the closed convex set  $X$ is
a hyperplane  which (i) has non-empty intersection with $X$,
and (ii) is the boundary of a closed half-space containing $X$. A {\it
  face} of $X$ is 
a subset $F$ of $\EE^n$ or $\HH^n$ 
which 
either is
  equal to $X$, 
or
is the intersection of $X$ with a support hyperplane for $X$. (In the
latter case we say that $F$ is a {\it proper} face of $X$.)
Note that a proper face of $X$ has strictly lower dimension than $X$.

We define a 
{\it convex \willbepolytope} in $\HH^n$ or $\EE^n$ to be
a  non-empty 
set which is 
the intersection of a family of closed
half-spaces whose bounding hyperplanes form a locally finite family.
Note that a face of a 
convex
\willbepolytope\ $X$ is itself a
convex
\willbepolytope\, and is
contained in $X$.

We define a {\it facet} of a convex polyhedron $X$ to be a maximal
proper face of $X$.


In this paper
we will use a number of 
facts about 
closed
convex sets and convex  \willbepolyhedra\ in $\HH^n$ that are the
counterparts of 
well-known facts about 
closed
convex  sets and convex \willbepolyhedra\  in
$\EE^n$.  
We will list several such facts here and then give very brief hints
about how to prove them.

\Claim\label{it's a ball}
Every $d$-dimensional closed
convex set 
$X$ in $\HH^n$ is homeomorphic to a set of the
form $D^d-\redXi $, where $D^d$ denotes the closed unit ball in $\RR^d$ and
$\redXi $ is some subset of $S^{d-1}=\partial D^d$. Furthermore, under a
homeomorphism between $X$ and $D^d-\redXi $,
the interior and boundary of $X$, defined as above, correspond respectively to $\inter
D^d$ and $S^{d-1}-\redXi $.
\EndClaim

\Claim\label{boundary}
The boundary of any 
closed
convex set is the union of its proper faces. 
\EndClaim

\Claim\label{good intersection}
If $X$ is a closed convex set, any 
non-empty set which is a
finite intersection of faces of $X$
is a face of $X$.
\EndClaim

\Claim\label{face of face}
If $X$ is a
convex \willbepolytope, any face of a face of $X$ is a face of $X$. 
\EndClaim 

\Claim\label{not-first fact} Every face of a convex \willbepolytope\ $X$ is a finite intersection of codimension-$1$ faces
of $X$. 
\EndClaim

\Claim\label{loc-fin faces}The faces of a convex \willbepolytope\ form
a locally finite family of sets. 
\EndClaim

Assertion \ref{it's a ball} 
is proved in the same way as its
Euclidean counterpart 
(which is a standard exercise).

Assertion \ref{good intersection}, which is the hyperbolic analogue of 
 \cite[Section 2.4, Assertion 10]{Grunbaum}, is easily proved
using the hyperboloid model of $\HH^n$. If $F_1,\ldots,F_m$ are faces
of the closed convex set $X\subset\HH^n$, we may write $F_i=\redPi_i\cap X$
for $i=1,\ldots,m$, where each $\redPi_i$ is a support hyperplane for
$X$. If we choose a point $x_0\in F_1\cap\cdots\cap F_m$, and use the
conventions of \cite[Section 1]{jason}, we may write each $\redPi_i$ as the
intersection of $\HH^n$ with a time-like hyperplane
$B_i=\{x\in\RR^{n+1}:u_i\circ(x-x_0)=0\}$, where $u_i\in\RR^{n+1}$ is a
vector such that $u_i\circ(x-x_0)\ge0$ for every $x\in X$. If we then
set $u=u_1+\cdots +u_m$, then $B:=\{x\in\RR^{n+1}:u\circ(x-x_0)=0\}$ is
time-like and $\redPi:=B\cap\HH^n$ is a support hyperplane for $X$, so that
$F:=F_1\cap\cdots\cap F_m=\redPi\cap X$ is a face of $X$.

In the special case of a convex polyhedron $X$ which is a {\it finite}
intersection of  closed half-spaces, the proofs of 
Assertions \ref{boundary}, 
\ref{face of face}, \ref{not-first fact}, and
\ref{loc-fin faces} above are almost identical to those of their
respective Euclidean analogues, \cite[Section 2.2, Assertion 4]{Grunbaum},
\cite[Section 2.6, Assertion 1]{Grunbaum}, \cite[Section 2.6,
Assertion 5]{Grunbaum}, and \cite[Section 2.6, Assertion
6]{Grunbaum}. The latter results concern closed
convex sets and so-called ``polyhedral sets'' in $\EE^n$; in
\cite{Grunbaum}, a ``polyhedral set'' is defined to be a subset of
$\EE^n$  which is a finite intersection of  closed half-spaces.

The
conclusions of the results cited from \cite{Grunbaum} are the same as
those of their hyperbolic analogues, except that the conclusion of \cite[Section 2.6, Assertion 6]{Grunbaum} gives
finiteness, rather than just local finiteness, of the set of faces. In
the proofs of the hyperbolic versions, one uses
Assertion \ref{good intersection} above in place of
 \cite[Section 2.4, Assertion 10]{Grunbaum}. Some trivial adjustments
 are needed, because the definitions given in \cite{Grunbaum} allow a
 ``polyhedral set'' or a ``face'' to be empty. Also, in adapting the
 proof of \cite[Section 2.6,
Assertion 5]{Grunbaum} to give a proof of Assertion \ref{not-first
  fact} above, one should 
replace the term``facet''  by ``codimension-$1$
  face''; the logic of the proof goes through equally well after this
  change. (We shall see in \ref{rebel} below that the notions of facet
  and of codimension-$1$ face are a posteriori equivalent.)

It is a straightforward exercise to extend the proofs of Assertions \ref{boundary}, 
\ref{face of face}, and \ref{loc-fin faces} to the  case of an arbitrary
convex polyhedron in $\HH^n$ from  the special case of a  finite
intersection of  closed half-spaces.
To extend 
\ref{not-first fact} to the general case, one needs the following
fact, which is a more interesting exercise. Let us say that 
a family $\redcalH$ of closed
half-spaces in $\HH^n$ is {\it irredundant} if $\redcalH$ has no proper subfamily $\redcalH'$ such that
$\bigcap_{\redrealcalH\in\redcalH'}\redrealcalH =\bigcap_{\redrealcalH\in\redcalH}\redrealcalH$.

\Claim\label{not so trivial}
If $\redcalH$ is a family of closed
half-spaces in $\HH^n$ whose bounding hyperplanes form a locally
finite family, then $\redcalH$ has an irredundant subfamily
$\redcalH_0$ such that
$\bigcap_{\redrealcalH\in\redcalH_0}\redrealcalH =\bigcap_{\redrealcalH\in\redcalH}\redrealcalH$.

In particular, every convex polyhedron in $\HH^n$ is the intersection
of an irredundant family of closed
half-spaces whose bounding hyperplanes form a locally
finite family.
\EndClaim

(In \cite{Grunbaum}, the case of the first assertion of \ref{not so
  trivial} in which the family $\redcalH$ is finite was left implicit, as
this case is trivial.)

The following facts will be useful.

\Claim\label{and this too}
The boundary of any convex \willbepolytope\ is the set-theoretical disjoint union of the interiors of all its proper faces.
\EndClaim

To prove \ref{and this too}, first note that an induction on the
dimension of a convex \willbepolytope\ $X$, using \ref{boundary} and \ref{face
  of face}, shows that $X$ is the union of the interiors of its
faces. To prove disjointness, suppose that $F_1$ and $F_2$ are
distinct faces
of $X$. By \ref{good intersection}, $G:=F_1\cap F_2$ is a face of $X$,
and the definition of a face then implies that $G$ is a face of $F_i$
for $i=1,2$. Since $F_1\ne F_2$, there is an index $i_0\in\{1,2\}$
such that $G$ is a proper face of $F_{i_0}$. We therefore have
$G\subset\partial F_{i_0}$, so that $\inter F_1\cap\inter
F_2=\emptyset$, and \ref{and this too} is proved.

\Claim\label{rebel}
A face of a convex polyhedron $X$ is a facet of $X$ if and only if it
has codimension $1$.
\EndClaim

The ``only if'' part of \ref{rebel} follows immediately from
\ref{not-first fact}. To prove the ``if'' part, suppose that $E$ is a
face of $X$  which is not a facet. Then $X$ has a proper
face $F$ with $E\subsetneq F$. This implies that $E$ is a proper face of $F$,
and hence $\dim E<\dim F<\dim X$; thus $E$ cannot have codimension $1$
in $X$.

\EndDefinitionsNotationRemarks

\DefinitionsNotationRemarks\label{complexes}

A {\it \willbepolytopal\ complex} in $\HH^n$ is a locally finite
collection $\realcalX$ of convex \willbepolyhedra\ such that (i) for each $X\in\realcalX$, all the faces of $X$ belong to $\realcalX$, and (ii) for all $X,Y\in\realcalX$, the set $X\cap Y$ either is a common face of $X$ and $Y$ or is empty. We denote by $|\realcalX|$ the union of all the convex \willbepolyhedra\ that belong to $\realcalX$.

The {\it open cells} of a \willbepolytopal\ complex $\realcalX$ are
the sets of the form $\inter X$ for $X\in\realcalX$. Note that 
by
\ref{it's a ball}, 
each
open cell of $\realcalX$ is topologically an open ball of some
dimension.

It follows from  \ref{and this too} and the definition of a
\willbepolyhedral\ complex 
that $|\realcalX|$ is set-theoretically the disjoint
union of the open cells of $\realcalX$.

\EndDefinitionsNotationRemarks

\DefinitionNotationRemarks\label{voronoi upstairs}
Let $\scrS$ be a locally finite (i.e. discrete and closed) subset of $\HH^n$. For each point $P\in \scrS$, we set
$$\XPS =\{W\in\HH^n:\dist(W,P)\le\dist(W,Q)\text{ for every
}Q\in \scrS\}.$$
We shall sometimes write $\calX_P$ in place of $\XPS $ in situations where
it is understood which  set $\scrS$ is involved.

A set having the form $\XPS $ for some $P\in \scrS$ will be called a
{\it Voronoi region} for the locally finite set $\scrS$.

\EndDefinitionNotationRemarks

\Proposition\label{first from jason}
Let $\scrS$ be a locally finite subset of $\HH^n$.
Then:
\begin{enumerate}
\item each Voronoi region
for $\scrS$ is a convex \willbepolytope;
\item if $X_0,\ldots,X_m$ are
Voronoi regions then $X_0\cap\cdots\cap X_m$ is a face of $X_i$ for
$i=0,\ldots,m$;
\item if $X$ is a
Voronoi region then every codimension-$1$ face of $X$ has the form
$X\cap Y$ for some Voronoi region $Y$; and
\item if $X_0$ is a
Voronoi region then every face of $X_0$ has the form
$X_0\cap X_1\cap\cdots\cap X_m$ for some Voronoi regions
$X_1,\ldots,X_m$.
\end{enumerate}
\EndProposition

\Proof
Assertion (1) is included in
\cite[Lemma 5.2]{jason}.
(The conventions of \cite{jason} regarding 
convex
polyhedra are the same
as those of the present paper, except that the term ``face''  is
defined in \cite{jason} to be what we call a proper face.)

Lemma 5.2 of \cite{jason} also asserts that if $m>0$, and if $X_0,\ldots,X_m$ are
distinct Voronoi regions for a locally finite set $\cals\subset\HH^n$, then
$X_0\cap\cdots\cap X_m$ is a proper face of $X_i$ for $i=0,\ldots,m$. This
implies (2).

To
prove (3), let $F$ be any codimension-$1$ face
of the Voronoi region 
$X$. 
Choose a point $U\in\inter F$. Since $U$
lies on the boundary of  $X$, it follows from the
definition of a Voronoi region (see \ref{voronoi upstairs}) that there is some Voronoi region $Y\ne
X$ such that $U\in Y$. Set $G=X\cap Y$; it follows from 
Assertion (2) of the present lemma
that $G$ is a face of $X$; it is a proper face since $Y\ne X$. By \ref{good intersection}, $L:=F\cap G$ is
a face of $X$; since $L$ is contained in the face $F$ of $X$, it
follows from the definition that $L$ is a face of $F$. But $L$
contains the interior point $U$ of $F$, and hence $L=F$. This means
that  $F\subset G$. 
But according to \ref{rebel},  the
codimension-$1$ face $F$ of $X$ is a facet, i.e. a maximal proper
face, of $X$; since $G$ is a proper face of $X$, 
it now follows that $F=G$, 
i.e. $F=X\cap Y$.

To  prove (4), suppose
that $E$ is a face of a Voronoi region $X$. 
According to
\ref{not-first fact}, 
we may write the face 
$E$ of $X$ as a finite intersection $F_1\cap\cdots\cap F_m$ of codimension-$1$ faces
of $X$. 
By Assertion (3), we may write $F_i=X\cap Y_i$ for $i=1,\ldots,m$, where
$Y_1,\ldots,Y_m$ are Voronoi regions. Then we have $E=X\cap
Y_1\cap\cdots\cap Y_m$, and the proof 
of (4) 
is complete.
\EndProof

\Notation
Let $\scrS$ be a locally finite subset of $\HH^n$. We denote by
$\realcalX_\scrS$ the set of all convex \willbepolyhedra\ which are faces of
Voronoi regions for $\scrS$.
\EndNotation

\Proposition\label{second from jason}
Let $\scrS$ be a locally finite subset of $\HH^n$. 
Then $\realcalX_\scrS$ is a
\willbepolytopal\ complex.
\EndProposition

\Proof
By 
Assertion (1) 
of Proposition \ref{first from jason}, each
Voronoi region is a convex \willbepolytope, and hence so are its faces. Thus
the elements of $\realcalX_\scrS$ are convex \willbepolyhedra.

To prove that the family $\realcalX_\scrS$ is locally finite, it
suffices---in view of \ref{loc-fin faces}---to prove that the Voronoi
regions
for $\scrS$ form a locally finite family. Thus for each
$Q\in\HH^n$, we must show that the set  $\scrT_Q$ of all
points $P\in\scrS$ such that  $\XPS\cap\nbhd_1(Q)\ne\emptyset
$ is finite. For this purpose, fix a point $P_0\in\scrS$. Given any point $P\in\scrT_Q$, we may
choose a point $R\in \XPS\cap\nbhd_1(Q)$. The definition of $\XPS$
gives $\dist(R,P)\le\dist(R,P_0)$. By the triangle inequality we then
have $\dist(Q,P)\le\dist(Q,P_0)+2\dist(R,Q)<\dist(Q,P_0)+2$. As the
latter inequality holds for every $P\in\scrT_Q$, the set $\scrT_Q$ is
bounded, and is therefore finite since $\scrS$ is locally finite. This
establishes the local finiteness of $\realcalX_\scrS$.

It follows from \ref{face of face} that any
face of a \willbepolytopein $\realcalX_\scrS$ is a \willbepolytopein
$\realcalX_\scrS$. Now suppose that $E$ and $E'$ are elements of
$\realcalX_\scrS$. By 
Assertion (4) 
of Proposition \ref{first from
  jason}, we may write $E=X_0\cap\cdots\cap X_m$ and
$E'=X_0'\cap\cdots\cap X_{m'}'$, where $X_0,\ldots, X_m$ and
$X_0',\ldots, X'_{m'}$ are Voronoi regions for $\scrS$.  According to 
Assertion (2)
of Proposition \ref{first from jason}, $E\cap
E'=X_0\cap\cdots\cap X_m\cap X_0'\cap\cdots\cap X_{m'}'$ is a face of
$X_0$. Since 
the face 
$E\cap E'$ of $X_0$
is contained in the face $E$ of $X_0$, it
follows from the definitions that $E\cap E'$
is a
face of $E$. Similarly it is a face of $E'$. This shows that 
$\realcalX_\scrS$  is a
\willbepolytopal\ complex.

\EndProof

\Remarks\label{bijection upstairs}
If $\scrS$ is a locally finite subset of $\HH^n$, then the 
$n$-dimensional open cells 
of $\realcalX_\scrS$ are precisely the interiors of the
Voronoi regions for $\scrS$.
It now follows from the definitions that for each $P\in \scrS$, the
open $n$-cell 
$\inter X_P$ is the unique open cell of $\realcalX_\scrS$ containing $P$, and that $P$ is the unique point of $\scrS$ lying in $\inter X_P$. Thus each point of $\scrS$ lies in a unique open $n$-cell of $\realcalX_\scrS$, and each open $n$-cell contains a unique point of $\scrS$.
\EndRemarks

\Definition\label{gpc def}
We define a {\it \gpc} to be an ordered quadruple
$K=(\redcalM,\cala,(\maybecalD_H)_{H\in\cala}, (\Phi_H)_{H\in\cala})$, 
where
\begin{enumerate}
\item $\redcalM$ is a Hausdorff space;
\item $\cala$ is a 
locally finite
collection of pairwise disjoint subsets of $\redcalM$
  whose union is $\redcalM$;
\item $(\maybecalD_H)_{H\in\cala}$ is an indexed family of convex \willbepolyhedra\,
  each contained in a hyperbolic space of dimension at least $2$;
\item $\Phi_H$ is a continuous map of $\maybecalD_H$ into $\redcalM$ for each
  $H\in\cala$;
\item $\Phi_H|\inter \maybecalD_H$ maps $\inter \maybecalD_H$ homeomorphically onto $H$
  for each $H\in\cala$; and
\item for each $H\in\cala$ and for each face $F$ of $\maybecalD_H$, there exist
  an element $C$ of $\cala$ and an
  isometry $\iota:F\to \maybecalD_C$ such that
  $\Phi_H|F=\Phi_C\circ\iota$. 
\end{enumerate}
Given a \gpc\ $K=(\redcalM,\cala,(\maybecalD_H)_{H\in\cala}, (\Phi_H)_{H\in\cala})$, we will call $\redcalM$ the {\it underlying space} of $K$ and will denote it
by $|K|$. The set $\cala$ will be denoted $\cala^K$, 
and its elements
 will be called {\it cells} of $K$. 
According to \ref{it's a ball}, each 
cell is
topologically an open ball of some dimension. For each $d\ge0$, the
set of all  $d$-cells, i.e. cells of dimension $d$, will be denoted by
$\cala^{d,K}$, and the union of all $d$-cells will be denoted by
$|\cala^{d,K}|$. For each cell $H$, the convex \willbepolytope\ $\maybecalD_H$ and
the map $\Phi_H$ will be denoted respectively by $\maybecalD_{H}^K$ and
$\Phi_{H}^K$. 
\invisiblecomment{
I've tried to fix subscripts and superscripts later to match this convention. I
think we're going to have dimensions always in superscripts, $K$
always a superscript, $S$ and $\tS$ always subscripts, $\realcalX$
always a superscript. In \ref{voronoi for M} we have defined
$\maybecalD_{H,S}$ to be $\maybecalD_H^{K_S}$, and have defined $\Phi_{H,S}$ similarly.
}

  Note that 
in the special case where $\maybecalD_H^K$ is compact for every $H\in\cala^K$, 
the
collection $\cala^K$ is in particular a CW complex with underlying space
$|K|$. 
Indeed, in this case,
according to \ref{it's a ball}, for each cell $H\in\cala^K$, 
the \willbepolytope\ $\maybecalD_H^K$ is topologically a closed ball whose
interior is mapped homeomorphically onto $H$ by $\Phi_H^K$; 
the map $\Phi_H^K$ therefore
gives rise to a characteristic map for $H$.
\EndDefinition

\Number\label{star}
If $K$
is a
\gpc, it
follows from Condition (6) of Definition \ref{gpc def} that for each
$H\in\cala^K$ and for each face $F$ of $\maybecalD_H^K$, the interior of $F$ is
mapped homeomorphically by $\Phi_H^K$ onto a cell of $K$. 
Given any cell
$C$ of $K$, we define the {\it star} of $C$ in $K$ to be the set of all
ordered pairs $(H,F)$ such that $H$ is a cell of $K$, $F$ is a face of
$\maybecalD_H^K$, and $\Phi_H^K(\inter F)=C$.
We define the {\it dimension} of an
element $(H,F)$ of the star of $C$ to be the dimension of $H$.

Note that if $K$ is a \gpc\ such that $|K|$ is an $n$-manifold without
boundary, the star of each $(n-1)$-cell of $K$ contains exactly two
elements of dimension $n$.
\EndNumber

\RemarksNotation\label{going down}
Let $n\ge2$ be an integer, and let $M$ be a   hyperbolic $n$-manifold.
Let us write $M=\HH^n/\Gamma$, where $\Gamma\le\isom(\HH^n)$ is
discrete and
torsion-free,
and let $q:\HH^n\to\Gamma$ denote the
quotient  map. 
Suppose that
$\realcalX$ is a
\willbepolytopal\ complex with $|\realcalX|=\HH^n$, such  that 
\begin{enumerate}
\item $\realcalX$ is
invariant under $\Gamma$ in the sense that $\gamma\cdot
X\in\realcalX$ whenever $\gamma\in\Gamma$ and
$X\in\realcalX$, and
\item
the stabilizer in $\Gamma$ of every element
of $\realcalX$ is trivial. 
\end{enumerate}
It follows from (1) and (2) 
that if $X$ is an arbitrary \willbepolytopein $\realcalX$, the restriction of $q$ to 
$\inter X$ is a homeomorphism of $\inter X$ onto a subset of $M$, which
will be denoted $C_X$; and that for 
any two elements $X,X'$ of
$\realcalX$ we have either $C_X=C_{X'}$ or $C_X\cap
C_{X'}=\emptyset$. Thus $\cala:=\{C_X:X\in\realcalX\}$ is a partition of
$M$;
it is locally finite since $q$ is a covering map, and since
$\realcalX$ is locally finite by the definition of a \willbepolytopal\
complex. 
For each $H\in\cala$, let us choose an element $\maybecalD_H$ of $\realcalX$
such that $C_{\maybecalD_H}=H$; and let us set $\Phi_H=q|\maybecalD_H$. Then 
$(M,\cala,(\maybecalD_H)_{H\in\cala}, (\Phi_H)_{H\in\cala})$ is a
\gpc\ with underlying space $M$. To verify Condition (6) of Definition
\ref{gpc def}, given an element $H$ of $\cala$ and a face $F$ of $\maybecalD_H$, we set $C=C_F$, so that $\maybecalD_C=\gamma(F)$ for some
$\gamma\in\Gamma$, and define the isometry $\iota$ to be
$\gamma|F$. The \gpc\ defined in this way will be denoted
$K^\realcalX$. It is well defined up to modifying the convex \willbepolyhedra\ $\maybecalD_H$ within
their isometry classes, and modifying the maps $\Phi_H$ by precomposing
them with isometries between convex \willbepolyhedra.
\EndRemarksNotation


\RemarksNotation\label{voronoi for M}
If $n\ge2$ 
is 
an integer, and  $S$ is a non-empty finite subset of
a  hyperbolic $n$-manifold $M$, we will associate with $S$ a \gpc\
$K_S$ such that $|K_S|=M$.
For this purpose we 
write $M=\HH^n/\Gamma$, where $\Gamma\le\isom(\HH^n)$ is discrete,
and torsion-free, and let $q:\HH^n\to
M
$ denote the
quotient  map. Then $\tS:=q^{-1}(S)$ is a non-empty, locally finite
subset of $\HH^n$. Hence by Proposition \ref{second
  from jason}, $\realcalX_{\tS}$ is a well-defined \willbepolytopal\
complex.
It is clear that the \willbepolytopal\ complex $\realcalX_{\tS}$ is invariant under
$\Gamma$, which is Condition (1) of \ref{going down} with $\realcalX=\realcalX_{\tS}$. We claim that
Condition (2) of \ref{going down} also holds with
$\realcalX=\realcalX_{\tS}$; to prove this, suppose that an element
$\gamma$ of $\Gamma$ stabilizes a polyhedron
$E\in\realcalX_{\tS}$. The set of Voronoi
regions
having $E$ as a face 
is non-empty by the definition of
$\realcalX_{\tS}$, and
is finite
since  $\realcalX_{\tS}$ is locally finite by the definition of a 
\willbepolytopal\ complex.
Since this set is invariant under $\gamma$, it
follows that there is an integer $m>0$ 
such that any Voronoi
region
$X$ 
having $E$ as a face is invariant under $\gamma^m$. But $\inter X$ contains a
unique point $P$ of $\tS$ by \ref{bijection upstairs}.
 Since $\tS$ is
$\Gamma$-invariant, it follows that $\gamma^m$ fixes $P$ and is
therefore the identity; as $\Gamma$ is torsion-free, this implies that
$\gamma=1$, and Condition (2) is established.

As Conditions (1) and (2) hold, \ref{going
  down} gives a well-defined  \gpc\ $K^{\realcalX_{\tS}}$ with
$|K^{\realcalX_{\tS}}|=M$. 
We now define $K_S$ to be
$K^{\realcalX_{\tS}}$. 

We shall also write
$\cala_S$ for
$\cala^{K_S}$,
and 
$\cala^d_S$ for
$\cala^{d,{K_S}}$ for any $d\ge0$. 
For each 
$H\in\cala_S$ we shall write
$\maybecalD_{H,S}$ for $\maybecalD_H^{K_S}$, and $\Phi_{H,S}$ for  $\Phi_H^{K_S}$.

According to \ref {bijection upstairs}, each point of $\tS$ lies in a
unique open $n$-cell of $\realcalX_{\tS}$, and each open $n$-cell
contains a unique point of $\tS$.
  Since, by \ref{going down},
every \willbepolytopein $\realcalX_{\tS}$ has trivial  stabilizer in $\Gamma$, it follows that 
each point of $S$ lies in a unique  $n$-cell of $K_S$, and that each $n$-cell of $K_S$ contains a unique point of $S$. 
In particular we have
$\#(\cala^n_S)=\#(S)$, and $S\subset\cala^n_S$.
\EndRemarksNotation

\Definition\label{dotsystem def}
Let $S$ be a non-empty finite subset of a  hyperbolic $3$-manifold. We
define a {\it \dotsystem} for $S$ 
to be a subset  $\redcalT $ of $|\atwoS|\subset K_S$ which contains at most one point
of each  $2$-cell of $K_S$.
\EndDefinition

\NotationRemarks\label{Gammadef}
Let $S$ be a non-empty finite subset of a  hyperbolic $3$-manifold,
and let $\redcalT $ be a 
\dotsystem\ for $S$. 
We will denote by $\redscrPST$  
the set of all
ordered pairs $(H,\redu)$ 
such that $H\in\athreeS$ and
$\redu\in{\Phi}_{H,S}^{-1}(\redcalT )\subset{\Phi}_{H,S}^{-1}(
|\atwoS|
)\subset\partial
\maybecalD_{H,S} 
$.
For each $H\in\athreeS$ we will denote by $\willbep_{H,\willbeS}$  the unique point of $S$
lying in $H$ (see \ref{voronoi for M}), and by $\willbehatp_{H,\willbeS}$  the
unique point of $\inter \maybecalD_{H,\willbeS}$ which is mapped to $\willbep_{H,\willbeS}$ by 
$\Phi_{H,S}$.

For
each 
$(H,\redu)\in\redscrPST$, let $\ell_{H,\redu}$ denote the line segment in the
convex \willbepolytope\ $ \maybecalD_{H,\willbeS}$ joining the points $\willbehatp_{H,\willbeS}$ and $\redu$. We
set
$$\scrG^{S,\redcalT }=S\cup\bigcup_{(H,\redu)\in
\redscrPST
}\Phi_{H,S}(\ell_{H,\redu})\subset
M.$$

If $\tau$ is a point in $\redcalT $, and if $C$ denotes the $2$-cell of $K_S$
containing $\tau$, then according to \ref{star}, the star of $C$ in $K_S$
has exactly two elements of dimension $3$. Hence there are exactly two
elements 
of $
\redscrPST
$, which we denote  $(H_i,\redu_i)$ for $i=1,2$, such that $\Phi_{H_i,S}(\redu_i)=\tau$. Now
$A_\tau:=
\{\tau\}
\cup\Phi_{H_1,S}(\inter\ell_{H_1,\redu_1})\cup\Phi_{H_2,S}(\inter\ell_{H_2,\redu_2})$
is an open arc in $M$. Topologically, $\scrG^{S,\redcalT }$ has the structure
of a graph (possibly with loops and multiple edges) in which the
vertices are the points of $S$ and the open edges are the arcs of the
form $A_\tau$ for $\tau\in \redcalT $.

It follows from this construction that for any $H\in\cala^3_S$, the
valence of the vertex 
$p_{H,\willbeS}$ 
in the graph $\scrG^{S,\redcalT }$ is
the 
number of
two-dimensional
faces of $ \maybecalD_{H,S}$ whose
interiors are mapped by $\Phi_{H,S}$ onto $2$-cells that meet 
$\redcalT $.
\EndNotationRemarks

\Proposition\label{same image}
Let $\Theta$ be 
a compact 
three-dimensional
submanifold-with-boundary of a 
hyperbolic $3$-manifold $M$. 
Let $S$ be a finite subset of 
$\Theta$, 
and let 
$\redcalT \subset M$ be a \dotsystem\ for $S\subset M$. Suppose that
for every  $2$-cell $C$ 
of $K_S$
with $C\cap\Theta\ne\emptyset$, we have $\redcalT \cap
C\cap\Theta\ne\emptyset$. 
Then the following conclusions hold.
\begin{enumerate}
\item  Any two points of $S$ that lie in the same component of
$\Theta$ also lie in the same component of the graph $\scrG^{S,\redcalT }$.
\item If $\willbep$ is a point of $S$, 
and if $\scrG_\willbep$ and $\Theta_\willbep$ respectively denote the
components of $\scrG^{S,\redcalT }$ and $\Theta$ that contain $\willbep$, then
$$\image(\pi_1(\Theta_\willbep,\willbep)\to\pi_1(M,\willbep))\le \image(\pi_1(\scrG_\willbep,\willbep)\to\pi_1(M,\willbep)),$$
where the unlabeled arrows denote inclusion
homomorphisms.
\end{enumerate}
\EndProposition

\Proof
Suppose that $M$, $\Theta$, $S$ and $\redcalT $ satisfy the hypotheses of the proposition.
We will prove the following assertion:

\Claim\label{ah forgot}
Let $\gamma:[0,1]\to\Theta$ be a  path in $\Theta$ whose endpoints lie
in $S$. Then $\gamma$ is fixed-endpoint homotopic in $M$ to a path in $\scrG^{S,\redcalT }$.
\EndClaim

Suppose for the moment that \ref{ah forgot} is true. If a point $p\in
S$ is given, applying \ref{ah forgot} to an arbitrary closed path in
$M$ based at $p$ establishes Conclusion (2) of the proposition. On the other hand, if $p$ and $q$ are points of $S$ that lie in the same component of
$\Theta$, applying   \ref{ah forgot} to a path in $\Theta$ joining $p$
to $q$ gives a path in $\scrG^{S,\redcalT }$ joining $p$
to $q$, and Conclusion (1) follows. Thus the proof of the proposition
will be complete once \ref{ah forgot} is proved.

Suppose that $\gamma$ satisfies the hypothesis of \ref{ah
  forgot}. Set $p_j=\gamma(j)$ for $j=0,1$.
Since $p_0,p_1\in S\subset|\athreeS|$ 
by \ref{voronoi for M}, we may assume after a
small fixed-endpoint homotopy that $\gamma$ maps $(0,1)$ into
$\inter\Theta$, that $\gamma([0,1])$ is disjoint from $|\anoughtS|$ and
$|\aoneS|$, and that $\gamma$ is transverse to 
the $2$-manifold
$|\atwoS|$.

If $\gamma^{-1}(|\atwoS|)=\emptyset$, then $\gamma([0,1])$ must be
contained in some $3$-cell $H$ of $K_S$. Since $H\cap S=\{p_{H,S}\}$ by \ref{Gammadef},
$\gamma$ is a closed path based at $p_{H,S}$, and in view of the
contractibility of $H$ it is fixed-endpoint homotopic to the constant path at
$p_{H,S}$. This gives the conclusion of \ref{ah forgot} in this case. 

For
the rest of the proof we will assume that $\gamma^{-1}(|\atwoS|)\ne\emptyset$.
We may then write 
$\gamma^{-1}(|\atwoS|)=\{s_i:0<i<n\}$, 
where
$n\ge2$
and $0=s_0<\cdots<s_n=1$. For $0<i<n$, let $C_i$ denote the
$2$-cell of $K_S$ containing $\gamma(s_i)$. Since $\gamma$ is a path
in $\Theta$, we have $C_i\cap\Theta\ne\emptyset$. The hypothesis of
the proposition then implies that $C_i$ contains a unique point $\redtau_i$
of $\redcalT $.

For each $i\in\{1,\ldots,n-1\}$ let us fix a path $\alpha_i$
in $C_i$ from $\gamma(s_i)$ to $\redtau_i$. 
We define  $\alpha_0$ and
$\alpha_n$ to be the constant paths at $p_0$ and $p_1$ respectively.

For $i=1,\ldots,n$ we have $\gamma((s_{i-1},s_i))\subset H_i$ for some
$H_i\in\athreeS$. 
Note that $p_0=p_{H_1,S}$ and that $p_1=p_{H_n,S}$.
We set $\maybecalD_i= \maybecalD_{H_i,\willbeS}$ and $\Phi_i=\Phi_{H_i,\willbeS}$ for
$i=1,\ldots,n$. 
There is a unique map  $\widehat\gamma_i$ from
$[s_{i-1},s_i]$ to the convex \willbepolytope\
$ \maybecalD_{\willbei}$
such that $\Phi_{\willbei}\circ\gamma_i=\gamma|[s_{i-1},s_i]$.

For $i=1,\ldots,n-1$ we have
$\gamma(s_{i})\in C_{i}$. Hence for $i=2,\ldots,n$ 
the point
$\widehat\gamma_i(s_{i-1})$ lies in the interior
of a  face
$\redU _i$
of $ \maybecalD_{\willbei}$, while for $i=1,\ldots,n-1$
the point
$\widehat\gamma_i(s_{i})$ lies in the interior
of a face
$\redU _i'$ of $ \maybecalD_{\willbei}$;
and
$\Phi_{\willbei}$  maps $\inter \redU _i$ 
homeomorphically onto $C_{i-1}$ for $i=2,\ldots,n$, while
$\Phi_{\willbei}$  maps $\inter \redU _i'$
homeomorphically onto $C_{i}$ for $i=1,\ldots,n-1$. 
There therefore exist
paths $\beta_i$ in  $\inter \redU
_i$
for
$i=2,\ldots,n$, and paths
$\beta_i'$ in $\inter \redU _i'$
for
$i=1,\ldots,n-1$, 
such
that $\Phi_{\willbei}\circ\beta_i=\alpha_{i-1}$ for $i=2,\ldots,n$ and
$\Phi_{\willbei}\circ\beta_i'=\alpha_{i}$ for $i=1,\ldots,n-1$. 
We define $\beta_1$ 
and $\beta_n'$ to be the constant paths at  $\hatp_0:=\hatp_{H_1,S}$ and $\hatp_1:=\hatp_{H_n,S}$
respectively; we then have
$\Phi_{\willbei}\circ\beta_i=\alpha_{i-1}$ and $\Phi_{\willbei}\circ\beta_i'=\alpha_{i}$ for $i=1,\ldots,n$.

For $i=1,\ldots,n$ we
may now choose a
homotopy $(\delta_{i,t}
:[s_{i-1},s_i]\to  \maybecalD_{\willbei}
)_{0\le t\le1}$ such that
$\delta_{i,0}=\hatgamma_i$, and such that
$\delta_{i,t}(s_{i-1})=\beta_i(t)$ and
$\delta_{i,t}(s_{i})=\beta_i'(t)$ for $0\le t\le1$. Now
$\epsilon_i:=\delta_{i,1}$ 
maps $s_{i-1}$ to $\redu_i:=\beta_i(1)$ and maps
$s_i$ to $\redu_i':=\beta'_i(1)$.

For $i=2,\ldots,n$ we have
$\Phi_{\willbei}(\redu_i)=\alpha_{i-1}(1)=\redtau_{i-1}\in \redcalT $; thus  we have
$(H_i,\redu_i)\in\redscrPST$, so that the line segment
$\ell_{H_i,\redu_i} \subset \maybecalD_{\willbei}$
is defined. Likewise, for $i=1,\ldots,n-1$
we have
$\Phi_{\willbei}(\redu_i')=\alpha_{i}(1)=\redtau_{i}\in \redcalT $; thus  we have
$(H_i,\redu_i')\in\redscrPST$, so that the line segment
$\ell_{H_i,\redu_i'}\subset \maybecalD_{\willbei}$ is defined. 
We set $L_1=\ell_{H_1,\redu_1'}$, $L_n=\ell_{H_n,\redu_n}$, and 
 $L_i=\ell_{H_i,\redu_i}\cup\ell_{H_i,\redu_i'}$ for any $i$ with
 $1<i<n$. Then $L_i\subset \maybecalD_{\willbei}$ is connected for
 $i=1,\ldots,n$; this is obvious for $i=1,n$, and for
$1<i<n$ it follows from the observation that
  the
segments $\ell_{H_i,\redu_i}$ and $\ell_{H_i,\redu_i'}$ share (at
least) the endpoint $\willbehatp_{H_i,\willbeS}$. (It may
happen that $\redu_i=\redu_i'$, but this does not affect the
argument.)

For $i=1,\ldots,n$, since $u_i$ and $u_i'$ lie in the connected set
$L_i$, we may choose a map $\zeta_i:[s_{i-1},s_i]\to L_i$ 
such that $\zeta_i(s_{i-1})=\redu_i$ and
  $\zeta_i(s_i)=\redu_i'$. Since $ \maybecalD_{\willbei}$ is contractible, there is a
  homotopy 
$(\eta_{i,t}
:[s_{i-1},s_i]\to  \maybecalD_{\willbei}
)_{0\le t\le1}$, constant on $\{s_{i-1},s_i\}$, such that
$\eta_{i,0}=\epsilon_i$ and $\eta_{i,1}=\zeta_i$.

For each $i\in\{1,\ldots,n\}$, let us now define $(\theta_{i,t}
:[s_{i-1},s_i]\to  \maybecalD_{\willbei}
)_{0\le t\le1}$ to be the composition of the homotopies $\delta_{i,t}$
and $\eta_{i,t}$ (so that $\theta_{i,t}=\delta_{i,2t}$ for $0\le
t\le1/2$ and $\theta_{i,t}=\eta_{i,2t-1}$ for $1/2\le
t\le1$). Then we have a well-defined homotopy $(\Lambda_{t}
:[0,1]\to M)_{0\le t\le1}$, constant on $\{0,1\}$, given by setting $\Lambda_t(s)=\Phi_{\willbei}\circ\theta_{i,t}(s)$
whenever $1\le i\le n$ and $s_{i-1}\le s\le s_i$. We have
$\Lambda_0=\gamma$, and the definition of $\scrG^{S,\redcalT }$ gives that
$\Lambda_1([0,1])\subset \scrG^{S,\redcalT }$.
\EndProof

\section{Volumes of certain convex sets in $\HH^3$} \label{convex section}

In the introduction to this paper we sketched the proof of the central
result Proposition \ref{central}. One essential step in the sketch involved computing the volume
of the set that was denoted in the introduction by $Z\cap 
\overline{\nbhd_\rho(\tredtau)}$.
Here $\rho$ is a positive real number and $\tredtau$ is a point in
$\HH^3$ (so that $\overline{\nbhd_\rho(\tredtau)}$ is a ball centered
at $\tredtau$); and $Z$ is an ``ice cream
cone'' defined as the convex hull of the union of $\{\tredtau\}$ with
another ball (a ``scoop of ice cream''). The goal of this section is to prove
Lemma \ref{Z-lemma}, which provides a formula for the volume of such
an intersection in what will turn out to be the crucial case: this is
the case in which the scoop has radius less than $\rho$, and the ball 
$\overline{\nbhd_\rho(\tredtau)}$ contains the center of the scoop
but does not contain the entire scoop. (These conditions are expressed
by a chain of inequalities in the hypothesis of Lemma \ref{Z-lemma}.)

Along the way to proving Lemma \ref{Z-lemma}, we will give formulae for the volumes
of various other kinds of sets in $\HH^3$.

\Notation\label{cap review}
We define  a 
strictly  increasing
function $B(\mayber)=\pi(\sinh(2\mayber)-2\mayber)$ for $\mayber>0$. Geometrically, $B(\mayber)$ gives the volume of a ball  in $\HH^3$ of
radius $\mayber$.

We define a {\it solid cap} to be a subset of $\HH^3$ which is the
intersection of a closed ball with a closed half-space which meets the
interior of the ball.
(In particular,  a closed ball is itself a solid cap in this
sense.)

Now let  real numbers $\mayber$ and $w$ be given, with $\mayber>0$.
Let $N$ be an open ball of radius $\mayber$ centered at a point 
$O\in \HH^3$, and let $\Pi$ be a plane whose distance from the center of
$N$ is $|w|$. 
Let $\redrealcalH$ be a closed half-space which is bounded by
$\Pi$, contains $O$ if $w\le0$, and does not contain $O$ if $w>0$. 
Then the set $\overline{N} \cap \redrealcalH$ is a 
solid cap if $|w|<\mayber$, and is otherwise either a one-point set or the empty set.
The volume of $\overline{N} \cap \redrealcalH$, which depends only on $\mayber$ and $w$, will be
denoted by $\kappa(\mayber,w)$.

For the case $w>0$, an analytic expression
for
$\kappa(r,w)$ is given in \cite[Section 14]{fourfree}. (Note that we have
$\kappa(\mayber,w)=0$ whenever $w\ge \mayber$.) We can calculate
$\kappa(r,w)$ in the case $w\le0$ by observing that
$\kappa(\mayber,0)=B(\mayber)/2$, and that
$\kappa(\mayber,w)=B(\mayber)-\kappa(\mayber,|w|)$ when $w<0$. This
method of calculating values of $\kappa$  is used in the course of
computations referred to in Section \ref{numerical section} of this
paper.
\EndNotation

Proposition \ref{altitude} below gives a formula for an altitude of a
hyperbolic triangle in terms of its sides; the function $\eta$ defined
in \ref{etadef} is needed for this formula.

\NotationRemarks\label{etadef}
We define a function 
$\eta$ 
on $(0,\infty)^3$ by
$$\eta(\wasx ,\wasy ,\wasz )=\frac{2(\cosh \wasx )(\cosh \wasy )(\cosh \wasz )-(\cosh^2\wasx  + \cosh^2 \wasy 
  +\cosh^2 \wasz )+1}{\sinh^2 \wasz }.
$$
(The labeling of the coordinates in $(0,\infty)^3$ as $\wasx$, $\wasy$
and $\wasz$ is chosen for consistency with the geometric application
to be given in Proposition \ref{altitude} below.)
Note that for any $\wasx ,\wasy ,\wasz 
>
0$ we have
$$(\sinh^2\wasz )((\sinh^2 \wasx )-\eta(\wasx ,\wasy ,\wasz ))=((\cosh \wasx )(\cosh \wasy )-\cosh \wasz )^2\ge0,$$
so that 
\Equation\label{mary moo}
\eta(\wasx ,\wasy ,\wasz )\le\sinh^2 \wasx 
\EndEquation
for all $\wasx ,\wasy ,\wasz >0$.

We set
$$\scrV=\{(\wasx ,\wasy ,\wasz )\in(0,\infty)^3:\eta(\wasx ,\wasy ,\wasz )\ge0\}.$$
Note that $\scrV$ is closed in the subspace topology of
$(0,\infty)^3\subset\RR^3$.

It follows from (\ref{mary moo}) that for every $(\wasx ,\wasy ,\wasz )\in\scrV$ we
have $(\cosh \wasx )/\sqrt{1+\eta(\wasx ,\wasy ,\wasz )}\ge1$. We may therefore define a non-negative-valued
function $\sigma$ on $\scrv$ by
$$\sigma(\wasx ,\wasy ,\wasz )=\arccosh\bigg(\frac{\cosh \wasx }
{\sqrt{1+\eta(\wasx ,\wasy ,\wasz )}}
\bigg).$$

\EndNotationRemarks


\Proposition\label{altitude}
Let $P_1P_2E$ be a hyperbolic triangle.
Set $r_i=|P_iE|$ for $i=1,2$, and $D=|P_1P_2|$. Let $\redU $ denote the perpendicular projection of $E$ to the line
containing $P_1$ and $P_2$ (so that $\redU $ may or may not lie on the
segment $P_1P_2$). 
Then 
$(r_1,r_2,D)\in\scrV$,
and $\sinh |E\redU |=\sqrt
{\eta(r_1,r_2,D)}$.
\EndProposition

\Proof
Let $\alpha$ denote the angle of the given triangle at the vertex
$P_1$. The hyperbolic law of cosines gives
\Equation\label{law of the jungle}
\cos\alpha=\frac{\cosh r_1\cosh D-\cosh r_2}{\sinh r_1\sinh 
D}.
\EndEquation
Writing $\sin^2\alpha=1-\cos^2\alpha$, substituting the right side of
(\ref{law of the jungle}) for $\cos\alpha$ and simplifying, we obtain
$\sin^2\alpha=\eta(r_1,r_2,D)/\sinh^2r_1$.
Hence
$\eta(r_1,r_2,D)\ge0$, i.e. $(r_1,r_2,D)\in\scrV$,
and 
$(\sinh r_1)(\sin\alpha)=\sqrt
{\eta(r_1,r_2,D)}$. But the hyperbolic law of sines, applied to the
right triangle $P_1\redU E$, gives $\sinh|E\redU |=\sinh r_1\sin\alpha$, and the
conclusion follows.
\EndProof

Our next goal is to give a formula for the intersection of two balls
in $\HH^3$.

\Notation\label{Vlens def}
We define a function $\Vlens$ on $\scrV$ by
$$\Vlens(x,y,z)=
\kappa(x,\sigma(x,y,z))+
\kappa(y,z-\sigma(x,y,z)).$$
\EndNotation

\Proposition\label{ball intersection}
Suppose that $P_1$ and $P_2$ are distinct points of $\HH^3$, and set
$D=\dist(P_1,P_2)$. Let positive numbers $r_1$ and $r_2$ be given, and
suppose that 
$r_2<\min(D,r_1)$, $D<r_1+r_2$, and $r_1<r_2+D$.
Then 
$(r_1,r_2,D)\in\scrV$, and
$$\vol(\overline{
\nbhd\nolimits_{r_1}(P_1)
}\cap \overline{
\nbhd\nolimits_{r_2}(P_2)
})=\Vlens(r_1,r_2,D).
$$
\EndProposition

\Proof
For $i=1,2$ set $N_i=\nbhd_{r_i}(P_i)$.
The hypothesis implies that
the open balls $N_1$ and $N_2$ have non-empty intersection and that
neither is contained in the other. Hence
$\Delta:=\partial\overline{N_1}\cap\partial\overline{N_2}$ is a closed disk. If 
$\Pi$
denotes the plane containing $\Delta $, the closed half-spaces bounded by
$\Pi$ may be labeled $\redrealcalH_1$ and $\redrealcalH_2$ 
in such a way that
$\overline{N_1}\cap\overline{N_2}$ is the union of the solid caps
$K_1:=\redrealcalH_1\cap \overline{N_1}$ and $K_2:=\redrealcalH_2\cap \overline{N_2}$.
We have $K_1\cap K_2=\Delta $. Hence
\Equation\label{told ya so}
\vol(\overline{N_1}\cap\overline{N_2})=\vol K_1+\vol K_2.
\EndEquation

Let $\redU $ denote the center of $\Delta $. 
Then  $\redU $, $P_1$ and $P_2$
lie on a line $L$, which meets 
$\Pi$ perpendicularly at $\redU $. We choose a point $E$ on the circle
$\partial \Delta $; then $\ell:=|E\redU |$ is the radius of $\Delta$.
For $i=1,2$, we denote by $u_i$ the distance from $P_1$
to the plane $\Pi$. 

We consider the right triangle $E\redU P_i$ for
$i=1,2$. Its hypotenuse has length
$|EP_i|=r_i$, and its other side lengths are $|E\redU |=\ell$ and
$|\redU P_i|=u_i$. 
The hyperbolic law of cosines therefore gives $\cosh
r_i=(\cosh u_i)(\cosh\ell)$, i.e. $r_i
>
\ell$ and
\Equation\label{still a thing}
u_i=\arccosh\bigg(\frac{\cosh r_i}{\cosh\ell}\bigg).
\EndEquation
The hypothesis implies that $r_2<r_1$, which with (\ref{still a
  thing}) gives
\Equation\label{horatio}
u_2<u_1.
\EndEquation

The triangle 
$P_1P_2E$ has side lengths
$r_1=|P_1E|$, $r_2=|P_2E|$ and $D=|P_1P_2|$; furthermore, the line
segment $E\redU $, whose length is $\ell$, is
contained in the plane $\Pi$ and is therefore perpendicular to $L$ at $\redU $.  It therefore follows from
Proposition \ref{altitude} that $(r_1,r_2,D)\in\scrV$ (which is the first
assertion of the present proposition), 
and that $\sinh \ell=\sqrt
{\eta(r_1,r_2,D)}$. Hence
\Equation\label{just a thing}
\cosh
 \ell=\sqrt
{1+\eta(r_1,r_2,D)}.
\EndEquation

It follows from (\ref{still a thing}), 
(\ref{just a thing}),
and the definition of the
function $\sigma$ on its domain $\scrV$
(see \ref{etadef}) that
\Equation\label{just when you thought}
u_1=\sigma(r_1,r_2,D).
\EndEquation

We claim that
\Equation\label{all day yesterday}
P_1\in\inter \redrealcalH_2=\HH^3-\redrealcalH_1.
\EndEquation
Assume that (\ref{all day yesterday}) is false, so that $P_1\in
\redrealcalH_1$. Since $P_1$ is the center of $N_1$, we then have $P_1\in
K_1\subset\overline{ N_1}\cap\overline{ N_2}$, so that $P_1\in\overline{ N_2}$. This means that
$\dist(P_1,P_2)\le r_2$, which is a contradiction, since by hypothesis
$\dist(P_1,P_2)=D>r_2$. Thus 
(\ref{all day yesterday}) is established.

Now since, by (\ref{all day yesterday}), 
the center of the radius-$r_1$ ball $N_1$, which lies at a 
distance $u_1$ from $\Pi$, does not lie in the half-space $\redrealcalH_1$,
it follows from the definition of
$\kappa$ given in \ref{cap review} that $\vol
K_1=\kappa(r_1,u_1)$. Combining this with 
(\ref{just when you thought}), we obtain
\Equation\label{poopageese}
\vol K_1=\kappa(r_1, \sigma(r_1,r_2,D)).
\EndEquation

Next note that since $\redU ,P_1,P_2\in L$ and $P_1\ne P_2$, one of the
following cases must occur: (a) $P_1$ lies strictly between $P_2$ and
$\redU $, or $\redU =P_1$; (b) $\redU $ lies strictly between $P_1$ and $P_2$; (c)
$P_2$ lies between $P_1$ and $\redU $, or $\redU =P_2$. If (a) holds, then
$u_2=\dist(\redU ,P_2)\ge\dist(\redU ,P_1)=u_1$, which contradicts
(\ref{horatio}). Hence (b) or (c) must hold. 

If (b) holds, then $P_1$ and $P_2$ lie in the interiors of distinct
half-spaces bounded by $\Pi$. Since $P_1\in\inter \redrealcalH_2$ by (\ref{all day
  yesterday}), we have $P_2\in\inter \redrealcalH_1$ in this case. 
Thus the center of the radius-$r_2$ ball $N_2$, which lies at a 
distance $u_2$ from $\Pi$, does not lie in the half-space $\redrealcalH_2$,
It therefore follows from the definition of
$\kappa$ given in \ref{cap review} that 
\Equation\label{ping-poing}
\vol
K_2=\kappa(r_2,u_2) \text{ if  (b) holds.}
\EndEquation

If (c) holds, then $P_1$ and $P_2$ lie in the same closed
half-space bounded by $\Pi$. Since $P_1\in\inter \redrealcalH_2$ by (\ref{all day
  yesterday}), we have $P_2\in \redrealcalH_2$ in this case. 
Thus the center of the radius-$r_2$ ball $N_2$, which lies at a 
distance $u_2$ from $\Pi$,  lies in the half-space $\redrealcalH_2$.
It therefore follows from the definition of
$\kappa$ given in \ref{cap review} that 
\Equation\label{table}
\vol
K_2=\kappa(r_2,-u_2) \text{ if (c) holds.}
\EndEquation

On the other hand, since $\dist(P_1,P_2)=D$, we have $u_2=D-u_1$ if
 (b) holds, and $u_2=u_1-D$ if  (c) holds. In view of (\ref{ping-poing})
and (\ref{table}), it follows that in all cases we have
$ \vol
K_2=\kappa(r_2,D-u_1)$.
Combining this with (\ref{just when you thought}),
we obtain
\Equation\label{i live in troy}
\vol K_2=\kappa(r_2,D-
\sigma(r_1,r_2,D))
\EndEquation
in all cases. The conclusion of the proposition follows from
(\ref{told ya so}), (\ref{poopageese}),
 (\ref{i live in troy}), and
the definition of $\Vlens$.
\EndProof

\Number\label{rcc}
By a {\it right circular cone} $\maybescrK\subset\HH^3$ we mean a set which is the
convex hull of $\{\maybeU\}\cup \Delta$, where $\maybeU$ (the {\it apex} of $\maybescrK$) is a
point, $\Delta$ is a closed disk contained in a plane $\Pi$, and the line
joining $\maybeU$ to the center of $\Delta$ is perpendicular to $\Pi$. A line segment joining $\maybeU$ to a point of $\partial \Delta$ will
be called a {\it generator} of $\maybescrK$; the line segment joining $\maybeU$ to the center of $\Delta$ will be called the {\it axis} of
$\Delta$; and the angle between a generator and the axis will be called the
{\it angle} of $\maybescrK$.
\EndNumber

The functions defined in the next subsection will be used in
Proposition \ref{rcc volume} to compute the volume of a right circular
cone in $\HH^3$.

\Notation
If $a$ is a positive number and $\beta$ is any real number, 
we set
$$\psi(a,\beta)=\arccosh\bigg(\frac{\cosh a}{\sqrt{1+(\sinh^2a)(\sin^2\beta)}}\bigg)$$
and
$$\Vcone(a,\beta)=\frac{B(a)}2(1-\cos\beta)-\kappa(a,\psi(a,\beta)).$$
\EndNotation

\Proposition\label{rcc volume}
If $\maybescrK$ is a right circular cone whose generator has length $a$ and
whose angle is $\beta$, then 
the axis of $\maybescrK$ has
length $\psi(a,\beta)$,
and we have 
$$\vol \maybescrK=\Vcone(a,\beta).$$
\EndProposition

\Proof
Let $\maybeU$ denote the apex of $\maybescrK$, let $O$ denote the center of
the base
$\Delta$ of $\maybescrK$, and let $P$ be a
point on $\partial \Delta$. Then $l:=|OP|$ is the radius of $\Delta$,
and $h:=|\maybeU O|$ is the length of the axis.
Applying the
hyperbolic law of sines to the right triangle $\maybeU OP$, we obtain
$\sinh l=(\sin\beta)(\sinh a)$; 
applying the hyperbolic Pythagorean theorem
to this same triangle, we obtain $(\cosh h)(\cosh l)=\cosh a$. In view
of the definition of $\psi(a,\beta)$ it follows that
$h=\psi(a,\beta)$, which is the first assertion of the proposition.

Now let $N$ denote the ball of radius $a$ centered at $\maybeU $, so that
$\vol N=B(a)$. Let $\Sigma$ denote the sector of $N$ consisting of all
radii of $N$ that meet $\Delta$. Then $\Sigma$ meets  the sphere
$S:=\partial N$ in a spherical cap centered at the point $Q$ where $S$
meets the
ray originating at $\maybeU $ and passing through $O$. Let us equip  $S$
with the standard spherical metric in which its total area is
$4\pi$. In terms of this metric we have
$S\cap\Sigma=\overline{\nbhd_\beta(Q)}$, so that the area of
  $S\cap\Sigma$ in the standard spherical metric is
  $2\pi(1-\cos\beta)$. Hence
$$\vol\Sigma=\bigg(\frac{\area(\Sigma\cap S)}{\area S}\bigg)\cdot B(a)=\frac{B(a)}2(1-\cos\beta).$$

Let $\redrealcalH$ and $\redrealcalH'$ denote the half-spaces in $\HH^3$ bounded by the
plane $\Pi$
containing $\Delta$, labeled in such a way that $O\in \redrealcalH'$. Since $h$ is the
distance from $\maybeU $ to $\Pi$, \ref{cap review} gives $\vol(\Sigma\cap
\redrealcalH)=\vol(N\cap \redrealcalH)=\kappa(a,h)$. But the sets $\Sigma\cap
\redrealcalH$ and $\maybescrK=\Sigma\cap
\redrealcalH'$ meet precisely in $\Delta$ and have union $\Sigma$; hence
$$\vol \maybescrK=\vol(\Sigma)-\vol(\Sigma\cap
\redrealcalH)=\frac{B(a)}2(1-\cos\beta)-\kappa(a,h);$$
since $h=\psi(a,\beta)$, the right-hand side of this last equality is
by definition equal to $\Vcone(a,\beta)$.
\EndProof

The next subsection introduces the functions that appear in Lemma
\ref{Z-lemma}, which is the main result of this section.

\Notation\label{phi-def}
If $r$ and $D$ are real numbers with $0<r<D$, 
we set
$$\omega(r,D)=\arccosh\bigg(\frac{\cosh
  D}{\cosh r}\bigg)$$
and
$$\theta(r,D)=\arcsin\bigg(\frac{\sinh r}{\sinh
  D}\bigg).$$
We then have $\omega(r,D)>0$ and $0<\theta(r,D)<\pi/2$.

We set
$\scrV_0=\{(\rho,r,D)\in\scrV:r<D\}$. Since $\Vlens$ has domain $\scrV$,
and since $\psi$, $\Vcone$ aand $\kappa$ all have domain
$(0,\infty)\times(-\infty,\infty)$, while $\omega(r,D)$ and $\theta(r,D)$ are
defined when $0<r<D$ and are positive-valued, we may define a function
$\phi$ with domain $\scrV_0$ by
$$\phi(\rho,r,D)=\Vlens(\rho,r,D)+\Vcone(\omega(r,D),\theta(r,D))
-\kappa(r,D-\psi(\omega(r,D),\theta(r,D))).$$
It follows from the definitions of the functions $\Vlens$, $\Vcone$,
$\kappa$ and $\psi$ that these functions are continuous on their
respective domains. Hence $\phi$ is continuous on $\scrV_0$.
\EndNotation

\Lemma\label{Z-lemma}
Let $\maybeU $ and $Q$ be distinct points in $\HH^3$. Set $D=\dist(\maybeU ,Q)$. 
Let
$\rho$ and $r$ be positive numbers such that
$r<D<\rho<D+r$.
Then $(\rho,r,D)\in\scrV_0$. Furthermore, if $Z\subset\HH^3$ denotes the convex hull
of $\{\maybeU \}\cup
\overline{\nbhd_{r}(Q)}
$, then 
$$\vol(Z\cap\overline{\nbhd\nolimits_\rho(\maybeU )})= 
\phi(\rho,r,D)
.$$
\EndLemma

\Proof
Set $J=\overline{
\nbhd_{r}(Q)}$.
Since $D>r$, there is a right
circular cone $\maybescrK$ which has apex $\maybeU $, and whose generators
are tangent to $J$.
We have $Z=J\cup\maybescrK$.
Let $a$ denote the common length of all generators
of $\maybescrK$.  If $\redU $ is a point of tangency
between a generator of
$\maybescrK$ and the ball $J$, then the triangle $\maybeU \redU Q$ has a right
angle at $\redU $; its hypotenuse $\maybeU Q$ has length $D$, and its sides $\maybeU \redU $
and $\redU Q$ have respective lengths $a$ and $r$. The hyperbolic
Pythagorean theorem therefore gives $\cosh D=(\cosh a)(\cosh r)$, i.e.
\Equation\label{mumbler}
a=\omega(r,D).
\EndEquation

In particular (\ref{mumbler}) implies that $a<D$, which in view of
the hypothesis implies that $a<\rho$. Since
$\maybescrK\subset\overline{\nbhd_a(\maybeU )}$, it follows that
$\maybescrK\subset\nbhd_{\rho}(\maybeU )$. Hence 
$$Z\cap\overline{\nbhd\nolimits_\rho(\maybeU )}=(J\cup\maybescrK)
\cap\overline{\nbhd\nolimits_\rho(\maybeU )}=
\maybescrK\cup(J\cap\overline{\nbhd\nolimits_\rho(\maybeU )}).$$
It follows that
\Equation\label{call me later}
\begin{aligned}
\vol(Z\cap\overline{\nbhd\nolimits_\rho(\maybeU )})&=\vol(\maybescrK)+\vol(J\cap\overline{\nbhd\nolimits_\rho(\maybeU )})-\vol(\maybescrK\cap(J\cap\overline{\nbhd\nolimits_\rho(\maybeU )}))\\
&=\vol(\maybescrK)+\vol(J\cap\overline{\nbhd\nolimits_\rho(\maybeU )})-\vol(\maybescrK\cap
J),
\end{aligned}
\EndEquation
where the last equality again follows from the inclusion $\maybescrK\subset\nbhd_{\rho}(\maybeU )$.

The angle $\beta$ of the right circular cone $\maybescrK$ is the angle
of the right triangle $\maybeU \redU Q$ at the vertex $\maybeU $. The hyperbolic law of
sines gives $\sin\beta=(\sinh r)/(\sinh D)$, i.e.
\Equation\label{tumbler}
\beta=\theta(r,D).
\EndEquation

According to Proposition \ref{rcc volume}, we have $\vol
\maybescrK=\Vcone(a,\beta)$, which with (\ref{mumbler}) and
(\ref{tumbler}) gives
\Equation\label{weird face}
\vol \maybescrK=\Vcone(\omega(r,D),\theta(r,D)).
\EndEquation
Since the hypotheses of the present lemma imply in particular that
$r<\min(D,\rho)$, $D<r+\rho$ and $\rho<r+D$, we may apply Proposition
\ref{ball intersection}, with $\maybeU $, $Q$, $\rho$ and $r$ playing the
respective roles of $P_1$, $P_2$, $r_1$ and $r_2$, to deduce that
$(\rho,r,D)\in\scrV$, and
that
\Equation\label{count of three}
\vol(J\cap\overline{\nbhd\nolimits_\rho(\maybeU )})=\Vlens(\rho,r,D).
\EndEquation

If $\Pi$ denotes the plane containing the base $\Delta$ of the cone
$\maybescrK$, it follows from the definition of $\maybescrK$ that
$\maybescrK\cap J$ is the intersection of $J$ with a half-space $\redrealcalH$
bounded by $\Pi$, and that the center $Q$ of the ball $J$ does not lie
in $\redrealcalH$. Hence if $d$ denotes the distance from $Q$ to $\Pi$, the
definition of $\kappa$ given in \ref{cap review} gives
$\vol(\maybescrK\cap J)=\kappa(r,d)$. But, again because $Q\notin \redrealcalH$,
the center $O$ of $\Delta$ lies on the line segment $Q\maybeU $, so that
$d=|QO|=D-|\maybeU O|$. Since $\maybeU O$ is the axis
of $\maybescrK$, it follows from the first
assertion of Proposition \ref{rcc volume} that
$|\maybeU O|=\psi(a,\beta)$. Hence
$\vol(\maybescrK\cap J)=\kappa(r,D-
\psi(a,\beta))$,
which with (\ref{mumbler}) and
(\ref{tumbler}) gives
\Equation\label{may be the last time}
\vol(\maybescrK\cap J)=\kappa(r,D-
\psi(\omega(r,D),\theta(r,D))).
\EndEquation

Since $r<D$ by hypothesis, and since we have observed that
$(\rho,r,D)\in\scrV$, we have
$(\rho,r,D)\in\scrV_0$. The equality
$\vol(Z\cap\overline{\nbhd\nolimits_\rho(\maybeU )})= 
\phi(\rho,r,D)
$
 follows from (\ref{call me later}),
(\ref{weird face}), (\ref{count of three}),
(\ref{may be the last time}), and the definition of $\phi(\rho,r,D)$.
\EndProof


\section{The central results}\label{central section}

The central results of this paper are Proposition \ref{central} and
its Corollary \ref {quote me}. The statements and proofs of these
results involve some functions that were used in \cite{boroczky}, of
which we review the definitions in Subsection \ref{what spoof}. They
also involve some basic notions about hyperbolic manifolds and metric
spaces  which are
discussed in Subsections \ref{nets and such one}---\ref{nets and such two}.

\NotationRemarks\label{what spoof}
As in \cite{boroczky}, for any $n\ge2$ and any
$r>0$ we shall denote by $h_n(r)$ the distance from the barycenter to
a vertex of a regular hyperbolic $n$-simplex $\scrD_{n,r}$ with sides
of length $2r$. 
Note that since $\scrD_{n,r}$ has diameter $2r$, we have $h_n(r)\le2r$
for any $n\ge2$ and any $r>0$.

Formulae for $h_2(r)$ and $h_3(r)$ are given in \cite[Subsection 9.1]{fourfree}.

For $r>0$
we define a function $\density(r)$ (denoted 
$d_3(r)$ in \cite{boroczky}),
by
$$\density(r)=(3\beta(r)-\pi)(\sinh(2r)-2r)/\tau(r),$$
where the functions 
$$\beta(r)=\arcsec(\sech(2r)+2)\text{ and }\tau(r)=3\int_{\beta(r)}^{\arcsec3}\arcsech((\sec t)-2)\,dt$$
respectively give the dihedral angle and the volume of $\scrD_{3,r}$.

For any $r>0$ we set $\maybeb(r)=B(r)/\density(r)$.

It is clear from the geometric definitions  that $\beta(r)$ and $\tau(r)$,
and therefore $\density(r)$ and $b(r)$, are continuous for $r>0$.

\EndNotationRemarks

\DefinitionRemark\label{nets and such one}
Let $X$ be a metric space, and let $\delta$ be a positive real
number. A {\it $\delta$-\magentapacking} for $X$ is defined to be a subset $S$ of
$X$ such that for any two distinct elements $x,y$ of $S$ we have
$\dist(x,y)\ge\delta$. A complete metric space $X$ is compact
if and only if for every $\delta>0$ there is a positive integer $N$
such that every $\delta$-\magentapacking\ in $X$ has cardinality at most $N$.
\EndDefinitionRemark

\DefinitionNotation\label{thin def}
A point $p$ of a hyperbolic $3$-manifold $M$ 
is said to be {\it $\maybeepsilon$-thin} for a given $\maybeepsilon>0$ if there is
a homotopically non-trivial loop based at $p$ having length less than
$\maybeepsilon$. The open set of $M$ consisting of all $\maybeepsilon$-thin points
is denoted $\Mthin(\maybeepsilon)$, and we set $\Mthick(\maybeepsilon)=M-\Mthin(\maybeepsilon)$.
\EndDefinitionNotation

\Definitions\label{marg def}
A {\it Margulis number} for an orientable hyperbolic $3$-manifold $M$
is defined to be a positive number $\tryepsilon$ such that, for every point
$p\in M$ and for any two loops $\alpha$ and $\beta$  based at $p$
such that the elements $[\alpha]$ and $[\beta]$ of $\pi_1(M,p)$ do not
commute, we have $\max(\length\alpha,\length\beta)\ge\tryepsilon$. If the
strict inequality $\max(\length\alpha,\length\beta)>\tryepsilon$ holds for all
$\alpha$ and $\beta$ such that $[\alpha]$ and $[\beta]$ do not
commute, we will say that $\tryepsilon$ is a {\it strict Margulis number} for
$M$.
\EndDefinitions

The following 
proposition
is a minor variant on a standard result:
see for example \cite[Proposition 4.9]{gs}. 
(The hypothesis that $\epsilon$ is a strict Margulis number, rather
than just a Margulis number, gives the additional information that $\Mthin(\tryepsilon)$ is the interior of a
smooth submanifold of $ M$.)
Most of the conclusions of this proposition
will be needed in the present paper; the final sentence will be needed
in one of the sequels to this paper.

\Proposition\label{strict summer vocation}
Let $M$ be a finite-volume orientable hyperbolic 
$3$-manifold, and let $\epsilon$ be a positive number. Then
$\Mthick(\epsilon)$ is compact. If in addition we assume that
$\tryepsilon$ is a strict Margulis number for $M$, then 
$\Mthin(\tryepsilon)$ is the interior of a
smooth submanifold of $ M$, which is closed as
a subset of $M$ and has only finitely many components; and each of these components is diffeomorphic to either $D^2\times S^1$ or
$T^2\times[0,\infty)$.
Hence $\Mthick(\tryepsilon)$ is a
connected 
$3$-manifold-with-boundary, 
and the inclusion homomorphism $\pi_1(\Mthick(\tryepsilon))\to\pi_1(M)$
is surjective. 
Furthermore, if $\tryepsilon$ is a strict Margulis number for $M$, and
if $M$ is written as a quotient $\HH^3/\Gamma$, where
$\Gamma\le\isomplus(\HH^3)$ is discrete and torsion-free with finite
covolume, and if $q:\HH^3\to M$ denotes the quotient map, then each
component of $q^{-1}(\Mthin(\epsilon))$ is convex.
\EndProposition

\Proof
Let $\delta>0$ be given. Set $\delta'=\min(\epsilon,\delta)$. If $S$
is any  $\delta$-\magentapacking\ in $\Mthick(\epsilon)\subset M$, and if $p_1,\ldots,p_m$ are
distinct points of $S$, then
the sets $\nbhd_{\delta'/2}(p_1),\ldots,\nbhd_{\delta'/2}(p_m)$ are
pairwise disjoint, and each of these sets
is isometric to a ball of radius
$\delta'/2$ in $\HH^3$. Hence $m\le\vol(M)/B(\delta'/2)$. This
shows that $\lceil\vol(M)/B(\delta'/2)\rceil$ is an upper bound for
the cardinality of any $\delta$-\magentapacking\ in $\Mthick(\epsilon)$, and hence $\Mthick(\epsilon)$ is compact.

To prove the second assertion, write $M=\HH^3/\Gamma$, where $\Gamma$ is a  torsion-free discrete
subgroup of $\isomplus(\HH^3)$, and let $q:\HH^3\to M$ denote the
quotient map. Suppose that $\tryepsilon$ is a strict Margulis number for $M$,
and for each non-trivial element $x$ of $\Gamma$ set
$Y(x)=\{P\in\HH^3:\dist(P,x\cdot P)
\le
\tryepsilon\}$. Let $\calc$ denote the set
of all maximal
abelian
subgroups  of $\Gamma$, and for each $C\in\calc$ set $W(C)=\bigcup_{1\ne x\in
  C}Y(x)$.  Now if
$C\in\calc$ is given, the non-trivial elements of $C$ all have 
either the
same axis $A_C$ or the same parabolic fixed point $T_C$;
and for
each $x\in C-\{1\}$, the set $Y(x)$ either is 
empty, or is  a  closed
hyperbolic cylinder centered at $A_C$ (i.e. a set of the form
$\{P\in\HH^3:\dist(P,A_C)\le R\}$ for some $R>0$), or is a closed
horoball based at $T_C$.
Furthermore, 
$C$ has only  finitely many elements $x$ such that
$Y(x)\ne\emptyset$. It follows that for each $C\in\calc$ we have
$W(C)=Y(x_C)$ for some $x_C\in C-\{1
\}
$. If $C$ and $C'$ are distinct
elements of $\calc$, then $x_C$ and $x_{C'}$ do not commute, and if
$P$ is  a point of $ {W(C)}\cap {W(C')}$, we have
$\max(\dist(P,x_C\cdot P), \dist(P,x_{C'}\cdot P))\le\tryepsilon$; this is
impossible since $\tryepsilon$ is a strict Margulis number for $M$. Hence
$(W(C))_{C\in\calc}$ is a disjoint family of subsets of $\HH^3$. It is
also a locally finite family, because if $C_1,C_2,\ldots$ is a
sequence of distinct elements of $\calc$, and $P_1,P_2,\ldots$ is a
sequence of points such that $P_i\in W(C_i)$ and $P_i\to P_\infty\in
\HH^3$, then for each $i\ge1$ there is a non-trivial element
$x_i$ of $C_i$ with $\dist(P_i,x_i\cdot P_i)\le\tryepsilon$; then
$x_1,x_2,\ldots$ is a sequence of distinct elements of $\Gamma$ for which
the sequence $(\dist(P_\infty,x_i\cdot P_\infty))_{i\ge1}$ is bounded,
contradicting the discreteness of $\Gamma$. Hence
$\bigcup_{C\in\calc}W(C)\subset\HH^3$ is a smooth 
$3$-manifold-with-boundary, 
closed as a subset of $\HH^3$, and its components are
cylinders and horoballs; it is clearly $\Gamma$-invariant, and
therefore has the form $q^{-1}(V)$ for some smooth manifold $V\subset
M$, which is closed as a subset of $M$. It follows from the definitions
that $\Mthin(\tryepsilon)=\inter V$. 
Each component of $V$ is identified with  $W(C)/C$ for some
$C\in\calc$. Since $W(C)/C$ is orientable and $\vol(W(C)/C)\le\vol
M<\infty$, the manifold $W(C)/C$ is diffeomorphic to 
$D^2\times S^1$  if $W(C)$ is a cylinder, and to
$T^2\times[0,\infty)$ if $W(C)$ is a horoball. The compactness of
$\Mthick(\tryepsilon)$ implies that $V$ has only finitely many components.

The final assertion now follows from the observation that a cylinder or
horoball in $\HH^3$ is convex.
\EndProof

\Definitions\label{nets and such two}
Let $M$ be a hyperbolic $3$-manifold and let $\maybeepsilon$ be a
positive real number. We define an {\it $\maybeepsilon$-thick $\maybeepsilon$-\magentapacking} for
$M$ to be an $\maybeepsilon$-\magentapacking\ for $M$ which is contained in
$\Mthick(\maybeepsilon)$. By a {\it maximal $\maybeepsilon$-thick $\maybeepsilon$-\magentapacking} for
$M$ we mean simply an $\maybeepsilon$-thick $\maybeepsilon$-\magentapacking\ for
$M$ which is not properly contained in any other $\maybeepsilon$-thick $\maybeepsilon$-\magentapacking\ for
$M$.
\EndDefinitions

In the statements and proofs of Proposition \ref{add this}, Lemma \ref{maximal properties} and
Proposition \ref{central}, we will freely use the machinery and
notation developed
in Section \ref{Voronoi section}: if
$\scrS$ is a locally finite subset of $\HH^3$, then for each point $P\in \scrS$, the set
$\XPS $ will be defined as in \ref{voronoi upstairs}, and will often
be denoted simply by $X_P$ when the set $\scrS$ is understood. If $S$ is a finite subset  of a 
hyperbolic $3$-manifold $M$, 
the \gpc\ $K_S$ will be defined as in
\ref{voronoi for M}, as will
the convex \willbepolytope\  $\maybecalD_{H,S}$  and the map $\Phi_{H,S}$
for
every cell $H$ of $K_S$, and the 
set $\cala^d_S$ of cells of $K_S$ for every integer $d\ge0$.

\Proposition\label{add this}
Let $\epsilon>0$ be given, and suppose that $\scrS$ is an
$\epsilon$-\magentapacking\ in $\HH^3$ (so that in particular $\scrS$ is a locally
finite subset of $\HH^3$). Then for every $P\in \scrS$, we have
$B(\epsilon/2)/\vol(\XPS\cap\nbhd_{h_3(\epsilon/2)}(P))\le
\density(\epsilon/2)$.
\EndProposition

\Proof
We shall extract this from the proof of
\cite[Theorem 1]{boroczky}. (The statement and proof of \cite[Theorem
1]{boroczky} are given in arbitrary spaces of constant curvature, but
in the following discussion we specialize to the hyperbolic case.)

If $\scrS$ is an $\epsilon$-\magentapacking\ in $\HH^n $, then for
each $P\in\scrS$, the ball $\nbhd_{\epsilon/2}(P)$ is contained in the
Voronoi region $\XPS$. In the terminology of \cite{boroczky}, the
balls $\nbhd_{\epsilon/2}(P)$ for $P\in\scrS$ are said to form a ``packing
of spheres of radius $\epsilon/2$.'' The set $\XPS$ is referred to as
the ``Voronoi-Dirichlet cell,'' or ``V-D cell,'' of the ball $\nbhd_{\epsilon/2}(P)$, and the ``density'' of
$\nbhd_{\epsilon/2}(P)$ in $\XPS$ is defined to be
$\vol(\nbhd_{\epsilon/2}(P))/\vol(\XPS)$. Theorem
1 of \cite{boroczky}
asserts that this ``density'' is bounded above by the quantity $d_n(\epsilon/2):=\vol(L\cap
\scrD_{n,\epsilon/2})/\vol(\scrD_{n,\epsilon/2})$, where $L$ denotes
the union of the balls of radius $\epsilon/2$ centered at the vertices
of $\scrD_{n,\epsilon/2}$. The quantity $d_3(\epsilon/2)$ is equal to
$\density(\epsilon/2)$ as defined in \ref{what spoof}.

In the proof of  \cite[Theorem 1]{boroczky}, $r$ denotes the quantity
which is denoted $\epsilon/2$ in the statement of the present lemma and
$S$ denotes the ball $\nbhd_{\epsilon/2}(P)$, while $\bar S$ denotes
the ball $\nbhd_{h_n(\epsilon/2)}(P)$ and $Z$ denotes the set $\XPS$.
The sentence in the proof beginning on line 3 of p. 259  reads:
``We shall show that the density of $S$ in $Z\cap\bar  S$ and a
fortiori in $Z$, is less than or equal to\ldots $d_n(r)$\ldots.'' In
our notation, the  ``density of $S$ in $Z\cap\bar  S$'' is
$\vol(\nbhd_{\epsilon/2}(P))/\vol(\XPS\cap\nbhd_{h_n(\epsilon/2)}(P))=B(\epsilon/2)
/\vol(\XPS\cap\nbhd_{h_n(\epsilon/2)}(P))$. 
Thus in the
three-dimensional case the latter quantity
is bounded above by $d_3(\epsilon/2)=\density(\epsilon/2)$.
\EndProof

\Lemma\label{maximal properties}
Let $M$ be a \willbefinitevolume\ 
orientable
hyperbolic $3$-manifold, and let $\maybeepsilon$ be a positive
number. Then: 
\begin{enumerate}
\item there exists a maximal $\maybeepsilon$-thick $\maybeepsilon$-\magentapacking\ for
$M$,
and any such 
set
is finite; and
\item if $S$ is any maximal $\maybeepsilon$-thick $\maybeepsilon$-\magentapacking\ for
$M$, then for every $3$-cell $H$ of $K_S$ we have 
$\vol H\ge
\maybeb(\maybeepsilon/2)$.
\end{enumerate}
 Furthermore, if we write $M=\HH^3/\Gamma$, where
 $\Gamma\le\isomplus(\HH^3)$  
is
 discrete
 and torsion-free
and has finite covolume, 
and let $q:\HH^3\to M$ denote the quotient map,
 then for any maximal $\maybeepsilon$-thick $\maybeepsilon$-\magentapacking\ $S$ for
$M$ we have
\begin{itemize}
\item[{\it(3)}]
$\tS:=q^{-1}(S)$ is an $\maybeepsilon$-\magentapacking\ for $\HH^3$,
\item[{\it(4)}] $\dist(\redE,\tS)<\maybeepsilon$ for every $\redE\in
  q^{-1}(\Mthick(\maybeepsilon))$; 
\item[{\it(5)}] for every $P\in\tS$ we have $\XPtS\supset\nbhd_{\maybeepsilon/2}(P)$ and $\XPtS\cap
  q^{-1}(\Mthick(\maybeepsilon))\subset\nbhd_{\maybeepsilon}(P)$;
and
\item[{\it(6)}] for every $P\in\tS$  we have
  $\vol(\XPtS\cap\nbhd_{h_3(\maybeepsilon/2)}(P))\ge
\maybeb(\maybeepsilon/2)
$.
\end{itemize}
\EndLemma

\Proof
Conclusion (1) follows from the compactness of $\Mthick(\epsilon)$,
which is the first assertion of Proposition \ref{strict summer vocation}.

To verify Conclusions 
(2)--(6), 
we take $M$ to be written in the form $\HH^3/\Gamma$, where $\Gamma$ is discrete
 and torsion-free, and we let $q:\HH^3\to M$ denote the quotient
 map. Set $\tS=q^{-1}(S)$. To prove (3) we must show that if  $P_1$ and $P_2$ are distinct points of
 $\tS$, then $\dist(P_1,P_2)\ge\maybeepsilon$. Set $p_i=q(P_i)$ for
 $i=1,2$. If $p_1\ne p_2$, then since $p_1$ and $p_2$ are distinct
 points of the $\maybeepsilon$-\magentapacking\ $S$, we have $\dist(p_1,p_2)\ge\maybeepsilon$ and
 hence $\dist(P_1,P_2)\ge\maybeepsilon$. If $p_1= p_2$, then  a geodesic path 
 from $P_1$ to $P_2$ is  projected by  $q$ to a closed path based at
 $p_1$, which is homotopically non-trivial since $P_1\ne P_2$. Since
 $p_1\in S\subset\Mthick(\maybeepsilon)$, this closed path must have length
 at least $\maybeepsilon$, and hence $\dist(P_1,P_2)\ge\maybeepsilon$ in this case
 as well.

To prove (4), let a point $\redE\in
  q^{-1}(\Mthick(\maybeepsilon))$ be given. Set $\rede=q(\redE)$. 
Since
  $\rede\in\Mthick(\maybeepsilon)$, we have $S\cup\{\rede\}\subset\Mthick(\maybeepsilon)$;
  the maximality of $S$ then implies that $S\cup\{\rede\}$ is not an
  $\maybeepsilon$-\magentapacking, and hence that $\dist(\rede,w)<\maybeepsilon$ for some $w\in
  S$. A path 
in $M$ from $\rede$ to $w$ which has 
length less than $\maybeepsilon$ lifts to a
  path of length less than $\maybeepsilon$ from $\redE$ to $W$ for some $W\in
  q^{-1}(\{w\})\subset q^{-1}(S)$, so that
  $\dist(\redE,\tS)\le\dist(\redE,W)<\maybeepsilon$, which gives (4).

To prove 
the first part of Assertion
(5), let $P\in\tS$ be given. Consider an arbitrary point
$\redE\in\nbhd_{\maybeepsilon/2}(P)$. If $P'$ is any point of $\tS$ distinct from
$P$, then by Conclusion (3) we have $\dist(P,P')\ge\maybeepsilon$. Hence
$\dist(\redE,P')\ge\dist(P,P')-\dist(\redE,P)>\maybeepsilon-\maybeepsilon/2=\maybeepsilon/2$, and
therefore $\dist(\redE,P)<\dist(\redE,P')$ for every point $P'\ne P$ in
$\tS$. 
By the definition of the Voronoi region $\willbeX_P$
 it follows
that $\redE\in \willbeX_P$; this establishes the inclusion
$\willbeX_P\supset\nbhd_{\maybeepsilon/2}(P)$.

To prove the second part of Assertion (5),
consider an arbitrary point $\redE\in \willbeX_P\cap
  q^{-1}(\Mthick(\maybeepsilon))$. Since $\redE\in 
  q^{-1}(\Mthick(\maybeepsilon))$, we have $\dist(\redE,\tS)<\maybeepsilon$ by
  Conclusion (4); that is, for some point $P'\in\tS$ we have
  $\dist(\redE,P')<\maybeepsilon$. Now since $\redE\in \willbeX_P$, we have
  $\dist(\redE,P)\le\dist(\redE,P')<\maybeepsilon$, so that
  $\redE\in\nbhd_{\maybeepsilon}(P)$. This establishes the inclusion
$\willbeX_P\cap
  q^{-1}(\Mthick(\maybeepsilon))\subset\nbhd_{\maybeepsilon}(P)$, and completes the
  proof of (5).

To prove (6), we apply Proposition \ref{add this}, letting $\tS$,
which is an $\epsilon$-\magentapacking\ by Assertion (3), play the role of $\scrS$
in that proposition. This gives $B(\epsilon/2)/\vol(\XPtS\cap\nbhd_{h_3(\epsilon/2)}(P))\le
\density(\epsilon/2)$, so that $\vol(\XPtS\cap\nbhd_{h_3(\epsilon/2)}(P))\ge
B(\epsilon/2)/\density(\epsilon/2)=b(\epsilon/2)$, and (6) is established.

Finally, to prove Conclusion (2), we note that according to the
definition of the \gpc\ $K_S$ (see \ref{going down} and \ref{voronoi for M}),
each cell $H$ of $K_S$ is the homeomorphic image under
$q$ of the Voronoi region $\willbeX_P$ for some $P\in\tS$. We then have
$\vol H=\vol \willbeX_P\ge \vol
(\willbeX_P\cap\nbhd_{h_3(\maybeepsilon/2)}(P))\ge
\maybeb(\maybeepsilon/2)
$, by Conclusion (6).
\EndProof

The following technical result, 
Lemma \ref{stoopid},
will be needed as background for the statement 
of
Proposition \ref{central}. 
The statements of  Lemma \ref{stoopid}, Proposition
\ref{central},  and Corollary \ref{quote me} involve the set $\scrV_0$ and the function $\phi$ that
were defined in Subsection \ref{phi-def}.

\Lemma\label{stoopid}
Let $\tryepsilon$ and $R$ be positive numbers such that $2\tryepsilon<R<5\tryepsilon/2$. Then for
every number $D$ in the interval 
$[R/2-\tryepsilon/4,\tryepsilon]$, 
we have
$(R-D,\tryepsilon/2,D)\in\scrV_0$. 
\EndLemma

\Proof
First consider the case of a point
$D\in(R/2-\tryepsilon/4,\tryepsilon]$. Set $\rho=R-D$ and $r=\tryepsilon/2$. The
inequalities $R>2\tryepsilon$ and
$R/2-\tryepsilon/4\le D\le 
\tryepsilon
$ then imply that
$r<D<\rho<D+r$.
It therefore follows from the first assertion of Lemma \ref{Z-lemma}
(applied to any two points of $\HH^3$ separated by a distance $D$)
that $(\rho,r,D)\in\scrV_0$. 

It remains to establish the assertion in the case
$D=D_0:=R/2-\tryepsilon/4$. We have shown that the continuous map $F:D\mapsto
(R-D,\tryepsilon/2,D)$ carries the interval $(R/2-\tryepsilon/4,\tryepsilon]$ into
$\scrV_0\subset\scrV\subset(0,\infty)^3$. 
The hypothesis of the present lemma directly implies that
$F(D_0)=(R/2+\tryepsilon/4,\tryepsilon/2,R/2-\tryepsilon/4)\in(0,\infty)^3$, and since $\scrV$
is closed in the subspace topology of $(0,\infty)^3\subset\RR^3$ (see \ref{etadef}), 
it
follows that $F(D_0)\in\scrV$. But the hypothesis also directly
implies that $\tryepsilon/2<R/2-\tryepsilon/4$, and so by definition
(see \ref{phi-def}) we have $F(D_0)\in\scrV_0$.
\EndProof

Among the conclusions of the following proposition, only Conclusion (5)
is quoted later in the paper (in the proof of Corollary \ref{quote me}). The other
conclusions are either needed for the proof of (5), or needed for
applications in the
sequel to this paper, or both.

\Proposition\label{central}
Let $M$ be a \willbefinitevolume\  orientable hyperbolic $3$-manifold. 
Let $\tryepsilon$ and
$R$ be positive numbers with 
$2\tryepsilon<R<5\tryepsilon/2$
(so that by Lemma \ref{stoopid}
we have 
$(R-D,\tryepsilon/2,D)\in\scrV_0$
for
every  
$D\in[R/2-\tryepsilon/4,\tryepsilon]$).
Let $c$ be a positive number such that (a) 
$\phi(R-D,
\tryepsilon/2,
D)\ge c$ for every $D\in[R/2-\tryepsilon/4,\tryepsilon]$, 
and (b) $B(\tryepsilon/2)\ge c$.

Set $\Theta=\Mthick(\tryepsilon)$, and let 
$S\subset\Theta$ be a maximal $\tryepsilon$-thick $\tryepsilon$-\magentapacking\ for $M$ (which
exists by Conclusion (1) of Lemma \ref{maximal properties}).
Let $\calc\subset\atwoS$ denote the set of all  $2$-cells of $K_S$
which meet $\Theta$. For each $C\in\calc$, select a point $\redtau_C\in
C\cap\Theta$, and 
define a \dotsystem\ $\redcalT $ for $S$ (see \ref{dotsystem def}) by setting
$\redcalT =\{\redtau_C:C\in\calc\}$.
Then:
\begin{enumerate}
\item The cardinality of $S$ is at most
$\lfloor\vol(M)/\maybeb(\tryepsilon/2)\rfloor$.
\item
For every $H\in\athreeS$, 
the number of
two-dimensional
faces of $ \maybecalD_{H,\willbeS}$ 
(see \ref{voronoi for M})
whose 
interiors are mapped by $\Phi_{H,S}$ onto cells in $\calc$
is at most 
$$\bigg \lfloor\frac{B(R)-\maybeb(\tryepsilon/2)}c
\bigg\rfloor.$$

\item
Each component of 
the graph $\scrG^{S,\redcalT }$  
(see \ref{Gammadef})
has first betti number bounded above
by
\Equation\label{gross}
1+ \bigg\lfloor\frac{\vol M}{\maybeb(\tryepsilon/2)}\bigg\rfloor\cdot\bigg(\frac12 \bigg \lfloor\frac{B(R)-\maybeb(\tryepsilon/2)}c
\bigg\rfloor
-1\bigg).
\EndEquation

\item
If $\tryepsilon$ is a 
strict
Margulis number for $M$, 
then 
$\scrG^{S,\redcalT }$ is connected, and
the inclusion homomorphism $\pi_1(\scrG^{S,\redcalT })\to\pi_1(M)$ is surjective.

\item
If  $\tryepsilon$ is a 
strict
Margulis number for $M$, we have
$$
\rank  \pi_1(M)\le1+ \bigg\lfloor\frac{\vol M}{\maybeb(\tryepsilon/2)}\bigg\rfloor\cdot\bigg(\frac12 \bigg \lfloor\frac{B(R)-\maybeb(\tryepsilon/2)}c
\bigg\rfloor
-1\bigg).
$$
\end{enumerate}
\EndProposition

\Proof
According to \ref{voronoi for M}, we have $\#(S)=\#(\athreeS)$. 
Since
the $3$-cells of $K_S$ are pairwise disjoint, we have
$\vol M\ge\sum_{H\in \athreeS}\vol H$. 
But by 
Conclusion (2) of Lemma \ref{maximal properties}
 we have 
$\vol H\ge
\maybeb(\tryepsilon/2)$
 for each $
H\in \athreeS$, and hence $\vol M\ge\#(
\athreeS
)\cdot
\maybeb(\tryepsilon/2)
=\#(
S
)\cdot \maybeb(\tryepsilon/2)
$. This implies (1).

To prove (2), 
we write $M=\HH^3/\Gamma$, where 
$\Gamma\le\isomplus(\HH^3)$ 
is discrete and
torsion-free,
and let $q:\HH^3\to M$ denote the
quotient  map. We set $\tS=q^{-1}(S)$.

If $H\in\cala_S^ 3$ is given, 
let us fix an open cell $\tH$ that is mapped homeomorphically onto $H$
by  $q$. We have $\tH=\inter \XPtS$ for some
$\willbeP\in\tS$.
The number of
two-dimensional
faces of $ \maybecalD_{H,\willbeS}$ 
whose 
interiors are mapped by $\Phi_{H,S}$ onto cells in $\calc$
is equal to $\#(\calf)$, where $\calf$ denotes the set of all 
two-dimensional
faces of
$\calX_{\willbeP}=\XPtS$
whose interiors are mapped by $q$ onto cells in $\calc$.
For each $F\in\calf$, 
we will denote by $T_F$
the unique point of $\inter F$ 
which is mapped to $\redtau_C$ by $q$.

Let
$F\in\calf$ be given. Then $F$
is a
codimension-$1$ face of the 
Voronoi region\ $\calX_P$. 
It
therefore follows from Assertion (3) of Proposition \ref{first from
  jason} that
 $F$ 
is the intersection of $\willbeX_{\willbeP}$ with another
Voronoi region,
which we write as $\willbeX_{Q_F}$ for
some $Q_F\in\tS$. 
(In the following argument, $X_P$ and $X_{Q_F}$ will play the role of the objects
that were denoted by $X$ and $Y$ when this argument was sketched in the introduction.)
In particular, $\inter\calX_{Q_F}$ is disjoint from
$\calX_{\willbeP}$ for each $F\in\calf$. On the other hand, if $F$ and $F'$
are distinct elements of $\calf$, the elements $\calX_{Q_F}$ and
$\calX_{Q_{F'}}$ of $\realcalX_{\tS}$ are distinct
and therefore have disjoint interiors. Thus $(\inter\calX_{Q_F})_{F\in\calf}$ is a disjoint family of subsets of $\HH^3-\calX_{\willbeP}$. 

Let $N\subset\HH^3$ denote the ball of radius $R$ centered at
$\willbeP$.
Set $\redL =N\cap X_P$.
By the hypothesis and an observation made in \ref{what spoof}, we have $R>2\tryepsilon>\tryepsilon\ge
h_3(\tryepsilon/2)$, 
so that $\nbhd_{h_3(\tryepsilon/2)}(P)\cap
X_P\subset \redL $. Since Conclusion (6) of Lemma \ref{maximal
properties} gives 
$\vol(\willbeX_P\cap\nbhd_{h_3(\tryepsilon/2)}(P))\ge
\maybeb(\tryepsilon/2)$, 
we have in particular that 
$\vol \redL \ge
\maybeb(\tryepsilon/2)$.
Hence
\Equation\label{first rat}
\vol(N-\redL )\le B(R)-\maybeb(\tryepsilon/2).
\EndEquation

Since $(\inter\calX_{Q_F})_{F\in\calf}$ is a disjoint family of
subsets of $\HH^3-\calX_{\willbeP}\subset\HH^3-\redL $, we have a disjoint
family $(N\cap\inter\calX_{Q_F})_{F\in\calf}$  of subsets of $N-\redL $. In
view of (\ref{first rat}) it follows that
\Equation\label{second rat}
\sum_{F\in\calf}\vol(N\cap\calX_{Q_F})\le B(R)-\maybeb(\tryepsilon/2).
\EndEquation

Now for each $F\in\calf$, 
the point
$T_F$ 
is a common point of the Voronoi regions
$\XPtS$ and 
$X_{Q_F,\tS}$, 
it follows from the definition
of a Voronoi
 region (see \ref{voronoi upstairs}) that
$\dist(\tredtau_F,\willbeP)=\dist(\tredtau_F,Q_F)$. We denote
the common value of these distances by $D_F$.

By hypothesis we have $\redtau_C\in C\cap\Theta$ for every $C\in\calc$. Hence
$\tredtau_F\in F\cap q^{-1}(\Theta)\subset X_P\cap q^{-1}(\Theta)$. But
according to Conclusion (5) of Lemma \ref{maximal properties} we have 
$X_P\cap
  q^{-1}(\Theta)\subset\nbhd_{\tryepsilon}(P)$. Thus
  $D_F=\dist(\tredtau_F,P)<\tryepsilon$. On the other hand, Conclusion (5) of Lemma
  \ref{maximal properties} also gives that 
$X_P\supset\nbhd_{\tryepsilon/2}(P)$, and since $\tredtau_F\in F\subset\partial
X_P$, it follows that $D_F=\dist(\tredtau_F,P)\ge\tryepsilon/2$. Thus we have
\Equation\label{ish}
\tryepsilon/2\le D_F\le\tryepsilon.
\EndEquation

The triangle inequality implies that 
\Equation\label{third rat}
\nbhd\nolimits_{R-D_F}(\tredtau_F)\subset\nbhd\nolimits_R(\willbeP)=N. 
\EndEquation

On the other hand, for each $F\in\calf$, we have $\tredtau_F\in
F\subset\calX_{Q_F}$, and  
Conclusion (5) of Lemma \ref{maximal properties}
implies that $
\overline{
\nbhd_{\tryepsilon/2}(Q_F)}
\subset\calX_{Q_F}$. Hence if we define 
$Z_F\subset\HH^3$ to be the convex hull of
$\{\tredtau_F\}\cup
\overline{
\nbhd_{\tryepsilon/2}(Q_F)}
$, the convexity of $\calX_{Q_F}$ (see
Proposition \ref{first from jason}, Assertion (1))
implies that $Z_F\subset\calX_{Q_F}$. With (\ref{third rat}), this gives
\Equation\label{fourth rat}
N\cap\calX_{Q_F}\supset
Z_F\cap \nbhd\nolimits_{R-D_F}(\tredtau_F).
\EndEquation

We claim that
\Equation\label{farfala}
\vol(Z_F\cap 
\clnbhd\nolimits_{R-D_F}(\tredtau_F))
\ge c.
\EndEquation
To prove (\ref{farfala}), we distinguish two cases. It follows from
(\ref{ish}) that either
$\tryepsilon/2\le D_F\le R/2-\tryepsilon/4$, or $R/2-\tryepsilon/4<D_F\le\tryepsilon$. To prove
(\ref{farfala}) in the case where $\tryepsilon/2\le D_F\le R/2-\tryepsilon/4$, we note
that in this case we have 
$R-D_F\ge D_F+\tryepsilon/2$; since $\dist(\tredtau_F,Q_F)=D_F$, it follows that $\nbhd_{R-D_F}(\tredtau_F)
\supset
\nbhd_{\tryepsilon/2}(Q_F)$. By definition the set $Z_F$ also contains
$\nbhd_{\tryepsilon/2}(Q_F)$, and hence $\vol(Z_F\cap \clnbhd_{R-D_F}(\tredtau_F))\ge
\vol \nbhd_{\tryepsilon/2}(Q_F)=B(\tryepsilon/2)$. But according to the hypothesis we
have $B(\tryepsilon/2)\ge c$, and (\ref{farfala}) is established in this case.

Now consider the case in which
$R/2-\tryepsilon/4<D_F\le\tryepsilon$.
If we set $D=D_F$, $\rho=R-D_F$ and $r=\tryepsilon/2$, the
inequalities $R>2\tryepsilon$ and
$R/2-\tryepsilon/4<D_F\le \tryepsilon$ then imply that 
$r<D<\rho<D+r$.
Since in addition we have $\dist(\tredtau_F,Q_F)=D$, 
we may apply Lemma \ref{Z-lemma}, with 
$Q_F$ and $T_F$  playing the respective roles of
 $Q$ and $T$ (so that $Z_F$ plays the role of $Z$) to deduce 
that
$(R-D_F,\tryepsilon/2,D_F)=(\rho,r,D)\in\scrV_0$, and that 
$\vol(Z_F\cap
\overline{\nbhd_{R-D_F}(\tredtau_F)})=\phi(\rho,r,D)=\phi(R-D_F,
\tryepsilon/2,
D_F)$. Since the boundary of the ball $\overline{\nbhd_{R-D_F}(\tredtau_F)}$
has volume $0$, we have $\vol(Z_F\cap
\nbhd_{R-D_F}(\tredtau_F))=\phi(R-D_F,
\tryepsilon/2,
D_F)$.
Since $D_F\in(R/2-\tryepsilon/4,\tryepsilon]\subset[R/2-\tryepsilon/4,\tryepsilon]$,Condition (a) of the hypothesis gives $\phi(R-D_F,
\tryepsilon/2,
D_F)\ge c$, and (\ref{farfala}) is established in this case as well. 

Now, combining 
(\ref{ish}) and (\ref{farfala}), we find
that
$\vol(N\cap\calX_{Q_F})\ge c$ 
for every $F\in\calf$. Hence the left hand side of (\ref{second rat}) is bounded below by $c\cdot\#(\calf)$. With (\ref{second rat}) we then obtain
$$\#(\calf)\le
\frac1c\cdot (B(R)-
\maybeb(\tryepsilon/2)),
$$
which gives (2).

To prove (3), let $\alpha$ denote the upper bound for $\#(S)$ given by
(1). 
According to \ref{Gammadef}, the vertex set of $\scrG^{S,\redcalT }$ is
equal to $S$, and hence has cardinality at most $\alpha$. Now let $\beta$ denote the quantity given by (2) which, for any
given 
$H\in\athreeS$,
bounds the 
number of
two-dimensional
faces of $ \maybecalD_{H,S}$ whose
interiors are mapped by $\Phi_{H,S}$ onto $2$-cells belonging 
to $\calc$.
Then according to \ref{Gammadef}, each vertex of $\scrG^{S,\redcalT }$ has
valence at most $\beta$.

In particular, if $\scrG$ is any component of $\scrG^{S,\redcalT }$ then the
set $\maybeL$ of vertices of $\scrG$ has cardinality at most $\alpha$,
and each  vertex of $\scrG$ has valence at most $\beta$. 
If $E$ denotes the number of edges of $\scrG$, the first betti number of $\scrG$ is 
$$1-\#(\maybeL)+E=1-\#(\maybeL)+\frac12\sum_{v\in \maybeL}{\rm valence}(v)=1+\sum_{v\in \maybeL}\bigg(\frac{{\rm valence}(v)}2-1\bigg)\le1+\alpha\bigg(\frac \beta2-1\bigg).$$
This proves (3).

Now assume that $\tryepsilon$ is a 
strict 
Margulis number for $M$. 
Then by 
Proposition \ref{strict summer vocation}, 
$\Theta$ is 
a connected 
$3$-manifold-with-boundary, 
and the inclusion homomorphism
$\pi_1(\Theta,\willbep)\to\pi_1(M,\willbep)$ is surjective. 
It follows from the construction of the \dotsystem\ $\redcalT $ in
the statement of the present proposition that for every  $2$-cell $C$ with $C\cap\Theta\ne\emptyset$
we have $\redcalT \cap C\cap\Theta\ne\emptyset$. Thus $M$, $\Theta$, $S$ and $\redcalT $
satisfy the hypotheses of 
Proposition
\ref{same image}. Since
$\Theta$ is connected, it follows from Assertion (1) of Proposition
\ref{same image} that all the points of $S$ lie in the same component
of $\scrG^{S,\redcalT }$. 
But according to 
\ref{Gammadef},
the vertex set of the graph $\scrG^{S,\redcalT }$ is $S$, and therefore every
component of $\scrG^{S,\redcalT }$ contains at least one point of $S$.  Hence $\scrG^{S,\redcalT }$ is connected.

Now let us choose
a point $\willbep\in S$.
Since $\Theta$ and $\scrG^{S,\redcalT }$ are connected,
Assertion (2) of Proposition \ref{same image} becomes 
$\image(\pi_1(\Theta,\willbep)\to\pi_1(M,\willbep))\le \image(\pi_1(\scrG^{S,\redcalT },\willbep)\to\pi_1(M,\willbep))$,
where the unlabeled arrows denote inclusion
homomorphisms.
Since the inclusion homomorphism
$\pi_1(\Theta,\willbep)\to\pi_1(M,\willbep)$ is surjective, it now
follows that the inclusion homomorphism
$\pi_1(\scrG^{S,\redcalT },\willbep)\to\pi_1(M,\willbep)$ is also surjective. 
This establishes (4).

To prove Assertion (5) we need only note that if $\tryepsilon$ is a strict
Margulis number for $M$, then it follows from
Assertion (4) that the rank of $\pi_1(M)$ is bounded above by the rank of
the fundamental group of the connected graph $\scrG^{S,\redcalT }$, which is in turn equal to the first betti number of
$\scrG^{S,\redcalT }$ 
and is therefore bounded above by (\ref{gross}) according to
Assertion (3).
\EndProof

\Corollary\label{quote me}
Let $M$ be a \willbefinitevolume\ orientable hyperbolic $3$-manifold. Let $\tryepsilon$ be
a (not necessarily strict) Margulis number for $M$, and let
$R$ be a number such that  $2\tryepsilon<R<5\tryepsilon/2$
(so that by Lemma \ref{stoopid}
we have 
$(R-D,\tryepsilon/2,D)\in\scrV_0$
for
every  $D\in[R/2-\tryepsilon/4,\tryepsilon]$).
Let $c$ be a positive number such that (a) 
$\phi(R-D,
\tryepsilon/2,
D)>c$ for every 
$D\in[R/2-\tryepsilon/4,\tryepsilon]$, 
(b) $B(\tryepsilon/2)> c$, and (c)
$(B(R)-\maybeb(\tryepsilon/2))/c$ is not an integer. Then we have
$$
\rank  \pi_1(M)\le1+ \bigg(\frac{\vol M}{\maybeb(\tryepsilon/2)}\bigg)\cdot\bigg(\frac12 \bigg \lfloor\frac{B(R)-\maybeb(\tryepsilon/2)}c
\bigg\rfloor
-1\bigg).
$$
\EndCorollary

\Proof
Choose a strictly
monotone increasing sequence $(\tryepsilon_i)_{i\ge1}$ converging to
$\tryepsilon$. Since $\tryepsilon$ is a Margulis number for $M$, and $\tryepsilon_i<\tryepsilon$ for
each $i\ge1$, each $\tryepsilon_i$ is a strict Margulis
number. Condition (b) of the hypothesis implies that $B(\tryepsilon_i/2)> c$
for all sufficiently large $i$. 

Since $R>2\tryepsilon_i$ for every $i$, it follows from Lemma \ref{stoopid} that
we have 
$(R-D,\tryepsilon_i/2,D)\in\scrV_0$
for
every $i\ge1$ and every $D\in[R/2-\tryepsilon/4,\tryepsilon]$.

We claim that for $i$ sufficiently
large, we also have $\phi(R-D,
\tryepsilon_i/2,
D)>c$ for every $D\in[R/2-\tryepsilon_i/4,\tryepsilon_i]$. If this is false, we may assume
after passing to a subsequence that for each $i\ge1$ there is a number
$D_i\in[R/2-\tryepsilon_i/4,\tryepsilon_i]\subset[R/2-\tryepsilon/4,\tryepsilon]$ such that
$\phi(R-D_i,
\tryepsilon_i/2,
D_i)\le c$. After again passing to a subsequence we may assume that the
sequence $(D_i)_{i\ge1}$ converges to a limit
$D_\infty\in[R/2-\tryepsilon/4,\tryepsilon]$. Lemma \ref{stoopid} gives
$(R-D_\infty,\tryepsilon/2,D_\infty)\in\scrV_0$, and since $\phi$ is
continuous on $\scrV_0$ (see \ref{phi-def}),
we have
$\phi(R-D_\infty,
\tryepsilon/2,
D_\infty)\le c$; this contradicts Condition (a) of the hypothesis, and
our claim is establshed.

Thus we have shown that for sufficiently large $i$, Conditions (a) and
(b) of the hypotheses
of Proposition \ref{central} hold when $\tryepsilon$ is replaced by
$\tryepsilon_i$. Since each $\tryepsilon_i$ is a strict Margulis number, it now
follows from Assertion (5) of Proposition \ref{central} that 
\Equation\label{take 'em}
\rank  \pi_1(M)\le1+ \bigg(\frac{\vol M}{\maybeb(\tryepsilon_i/2)}\bigg)\cdot\bigg(\frac12 \bigg \lfloor\frac{B(R)-\maybeb(\tryepsilon_i/2)}c
\bigg\rfloor
-1\bigg)
\EndEquation
for sufficiently large $i$. But since the function $\maybeb$ is
continuous
by \ref{what spoof},
we have $(B(R)-b(\tryepsilon_i/2))/c\to (B(R)-b(\tryepsilon/2))/c$ as $i\to\infty$;
and Condition (c) of the hypothesis then guarantees that
$\lfloor(B(R)-b(\tryepsilon_i/2))/c\rfloor=
\lfloor(B(R)-b(\tryepsilon/2))/c\rfloor$ for sufficiently large $i$. The conclusion of the corollary
therefore follows upon taking limits in (\ref{take 'em}).
\EndProof

\section{Numerical calculations}
\label{numerical section}

In this section we apply Corollary \ref{quote me} to the proofs of
Proposition \ref{pi-one rank} and Theorem \ref{main}, which provide
the concrete estimates stated in the introduction.

The statement and proof of the following lemma
involve the set $\scrV_0$, and the functions
$\eta$, 
$\sigma$, 
$\Vlens$, $\omega$, $\theta$,
$\psi$, $\Vcone$
and
$\phi$, which  were
defined in Section \ref{convex section}.

\Lemma\label{compute}
Set $\tryepsilon=\log3$ and $R=2\log3+0.15$. Let $I$ denote the interval
$[R/2-\tryepsilon/4,\tryepsilon]=[3(\log3)/4+.075,\log3]$, so that by Lemma
\ref{stoopid} we have
$(R-D,\tryepsilon/2,D)\in\scrV_0$ for every $D\in I$. Then we have $\phi(R-D,\tryepsilon/2,D)>0.496$ for every $D\in I$.
\EndLemma

\Proof
Since $(R-D,\tryepsilon/2,D)\in\scrV_0$ for every $D\in I$, we may define
functions $H$, $\Sigma$, $\Wlens$, $\Psi$, 
$\Wcone$, and $\Phi$ on $I$ by $H(D)=\eta(R-D,\tryepsilon/2,D)$, 
$\Sigma(D)=\sigma(R-D,\tryepsilon/2,D)$, 
$\Wlens(D)=\Vlens(R-D,\tryepsilon/2,D)$,
$\Psi(D)=\psi( \omega(\tryepsilon/2,D), \theta(\tryepsilon/2,D))$,
$\Wcone(D)=\Vcone( \omega(\tryepsilon/2,D), \theta(\tryepsilon/2,D))$,
and
$\Phi(D)=\phi(R-D,\tryepsilon/2,D)$. We shall begin by finding upper and lower bounds for 
$H$, $\Sigma$,  and $\Psi$, and lower bounds for
$\Wcone$, 
$\Wlens$, and $\Phi$,
on certain subintervals of $I$.

From the definitions of $H$, $\Sigma$, $\Wlens$, $\Psi$, 
$\Wcone$, and $\Phi$ given above, and the definitions of $\eta$,
$\sigma$, $\Vlens$, $\psi$, $\Vcone$ and $\phi$ given in Section
\ref{convex section}, we find that, for each $D\in I$, we have

\Equation\label{Hi I'm H}
H(D)=\frac{2\cosh (R-D)\cosh(\tryepsilon/2)\cosh D-(\cosh^2(R-D) + \cosh^2(\tryepsilon/2)
  +\cosh^2 D)+1}{\sinh^2 D},
\EndEquation

\Equation\label{Hi I'm Sigma}
\Sigma(D)=
\arccosh\bigg(\frac{\cosh(R-D)}
{\sqrt{1+H(D)}}.
\bigg),
\EndEquation

\Equation\label{Hi I'm Wlens}
\Wlens(D)=
\kappa(R-D,\Sigma(D))+
\kappa(\tryepsilon/2,D-\Sigma(D)),
\EndEquation

\Equation\label{Hi I'm Psi}
\Psi(D)=\arccosh\bigg(\frac{\cosh\omega(\tryepsilon/2,D)
}{\sqrt{1+(\sinh^2\omega(\tryepsilon/2,D))(\sin^2\theta
(\tryepsilon/2,D)
)}}\bigg),
\EndEquation

\Equation\label{Hi I'm Wcone}
\Wcone(D)=\frac{B(\omega(\tryepsilon/2,D)
)}2(1-\cos
\theta(\tryepsilon/2,D)
)-\kappa(
\omega(\tryepsilon/2,D),
\Psi(D)
)),
\EndEquation

and

\Equation\label{Hi I'm Phi} 
\Phi(D)=\Wlens(D)+\Wcone(D)
-\kappa(\tryepsilon/2,D-\Psi(D)).
\EndEquation

For any nondegenerate closed subinterval $[D_-,D_+]$ of $I$, we set 
$$H^-(D_-,D_+)=\frac{2\cosh (R-D_+)\cosh(\tryepsilon/2)\cosh D_--(\cosh^2(R-D_-) + \cosh^2(\tryepsilon/2)
  +\cosh^2 D_+)+1}{\sinh^2 D_+}
$$
and
$$H^+(D_-,D_+)=\frac{2\cosh (R-D_-)\cosh(\tryepsilon/2)\cosh D_+-(\cosh^2(R-D_+) + \cosh^2(\tryepsilon/2)
  +\cosh^2 D_-)+1}{\sinh^2 D_-}.
$$
These definitions, together with (\ref{Hi I'm H}), imply that
\Equation\label{H is good}
H^-(D_-,D_+)\le  H(D) \le H^+(D_-,D_+) \text{ for every }D\in[D_-,D_+].
\EndEquation

We define a {\it good interval} to be a nondegenerate closed
subinterval $[D_-,D_+]$ of $I$ such that (1) $H^-(D_-,D_+)>-1$, (2)
$H^+(D_-,D_+)<
\sinh^2(R-D_+)$ and (3) $\sinh \omega(\tryepsilon/2,D_-)>\sinh
\omega(\tryepsilon/2,D_+)\cdot\sin \theta(\tryepsilon/2,D_+)$.

Let $[D_-,D_+]$ be any good interval. It follows from Conditions (1)
and (2) of the definition of a good interval, together with (\ref{H is
  good}), that $-1<H^-(D_-,D_+)\le
H^+(D_-,D_+)<\sinh^2(R-D_+)<\sinh^2(R-D_-)$. Hence we may define quantities
$\Sigma^-(D_-,D_+)$ and $\Sigma^+(D_-,D_+)$ by setting
$$\Sigma^-(D_-,D_+)=
\arccosh\bigg(\frac{\cosh(R-D_+)}
{\sqrt{1+H^+(D_-,D_+)}}
\bigg)
$$
 and 
$$\Sigma^+(D_-,D_+)=
\arccosh\bigg(\frac{\cosh(R-D_-)}
{\sqrt{1+H^-(D_-,D_+)}}
\bigg).
$$
These definitions, together with (\ref{Hi I'm Sigma}) and (\ref{H is
    good}), imply
that 
\Equation\label{the real mcshnuff}
\Sigma^-(D_-,D_+)\le\Sigma(D)\le\Sigma^+(D_-,D_+)
\EndEquation
for every good interval $[D_-,D_+]$ and every $D\in [D_-,D_+]$.

Now, for every good interval $[D_-,D_+]$,  set 
$$
\Wlens^-(D_-,D_+)=
\kappa(R-D_+,\Sigma^+(D_-,D_+))+
\kappa(\tryepsilon/2,D_+-\Sigma^-(D_-,D_+)).$$
Since
the function $\kappa$ is (weakly)
monotone
increasing in its first (positive-valued) argument and monotone
decreasing in its second (real-valued) argument, it follows from (\ref{Hi I'm Wlens}),
(\ref{the real mcshnuff}), and the definition of $\Wlens^+$ that for every good interval
$[D_-,D_+]$ and every $D\in[D_-,D_+]$ we have
\Equation\label{magootz}
\Wlens(D)\ge\Wlens^-(D_-,D_+).
\EndEquation

Recall from  \ref{phi-def} that the function $\theta$ takes its values in $(0,\pi/2)$.
 Condition (3) of the definition of a good interval, together with the
 positivity of the sine function on $(0,\pi/2)$,
implies 
that for any good interval $[D_-,D_+]$ we have
$$
\begin{aligned}
1+(\sinh^2
\omega(\tryepsilon/2,D_-))(\sin^2 \theta(\tryepsilon/2,D_-))&<
1+(\sinh^2
\omega(\tryepsilon/2,D_+))(\sin^2 \theta(\tryepsilon/2,D_+))\\
&<\cosh^2
\omega(\tryepsilon/2,D_-)
<\cosh^2
\omega(\tryepsilon/2,D_+).
\end{aligned}
$$
Hence for each good interval
$[D_-,D_+]$ we may define  quantities $\Psi^-(D_-,D_+)$ and
$\Psi^+(D_-,D_+)$ by
$$\Psi^-(D_-,D_+)=\arccosh\bigg(\frac{\cosh\omega(\tryepsilon/2,D_-)
}{\sqrt{1+(\sinh^2\omega(\tryepsilon/2,D_+))(\sin^2\theta
(\tryepsilon/2,D_+)
)}}\bigg)
$$
and
$$\Psi^+(D_-,D_+)=\arccosh\bigg(\frac{\cosh\omega(\tryepsilon/2,D_+)
}{\sqrt{1+(\sinh^2\omega(\tryepsilon/2,D_-))(\sin^2\theta
(\tryepsilon/2,D_-)
)}}\bigg).
$$
This definition, together with (\ref{Hi I'm Psi}), the monotonicity
of the sine function on $(0,\pi/2)$, and the observation that $\omega$
and $\theta$ (defined in \ref{phi-def}) are monotone increasing in
their second argument,  implies that 
\Equation\label{Psi is good}
\Psi^-(D_-,D_+)\le  \Psi(D) \le \Psi^+(D_-,D_+)
\EndEquation
for every good interval  $[D_-,D_+]$ and   for every $D\in[D_-,D_+]$.

Next, for every good interval $[D_-,D_+]$, set 
$$\Wcone^-(D_-,D_+)=\frac{B(\omega(\tryepsilon/2,D_-)
)}2(1-\cos
\theta(\tryepsilon/2,D_-)
)-\kappa(
\omega(\tryepsilon/2,D_+),
\Psi^-(D_-,D_+)
).
$$
It follows from this definition, together with (\ref{Hi I'm Wcone}), (\ref{Psi is good}),
the monotonicity properties of $\kappa$, $\omega$ and $\theta$
mentioned above, and the monotonicity of the cosine function on
$(0,\pi/2)$, that 
\Equation\label{Wcone is good}
\Wcone(D)\ge\Wcone^-(D_-,D_+)
\EndEquation
for every good interval  $[D_-,D_+]$ and   for every $D\in[D_-,D_+]$.

Now, for every good interval $[D_-,D_+]$, set
$$\Phi^-(D_-,D_+)=
\Wlens^-(D_-,D_+)+\Wcone^-(D_-,D_+)
-\kappa(\tryepsilon/2,D_- -\Psi^+(D_-,D_+)).$$

It follows from this definition, together with
(\ref{magootz}), (\ref{Wcone is good}), (\ref{Psi is good}), and the monotonicity properties of $\kappa$
mentioned above, that 
\Equation\label{Phi is good}
\Phi(D) \ge \Phi^-(D_-,D_+)
\EndEquation
for every good interval  $[D_-,D_+]$ and   for every $D\in[D_-,D_+]$.

In view of (\ref{Phi is good}), in order to prove the lemma it
suffices to show that every point of $I$ lies in a good interval
$[D_-,D_+]$ such that $\Phi^-(D_-,D_+)>0.496$.

For this purpose, we denote by $(\delta_1,\ldots,\delta_{46})$ the
$46$-tuple of positive constants
$$
\begin{aligned}
(&0.17,0.14,0.12,0.10,0.09,0.08,0.07,0.06,0.05,0.045,0.040,0.035,0.030,0.025,\\&0.022,0.020,0.018,0.016,0.014,0.012,0.010,0.0084,0.007,0.006,0.005,0.0042,\\&0.0035,0.0030,0.0025,0.0022,0.0019,0.0016,0.0013,0.0011,0.0009,0.00075,\\&0.0006,0.0005,0.0004,0.0003,0.00025,0.00020,0.00015,0.00010,0.00005,0.00002).
\end{aligned}
$$
We set $D_0=R/2-\tryepsilon/4$, $D_i=\tryepsilon-\delta_i$ for $i=1,\ldots,46$ and
$D_{47}=\tryepsilon$. We have $R/2-\tryepsilon/4=D_0<D_1<\cdots<D_{47}=\tryepsilon$, so that $I$
is the union of the intervals $[D_0,D_1],\cdots,[D_{46},D_{47}]$. By direct
calculation we find that for $i=1,\ldots,47$ the interval
$[D_{i-1},D_{i}]$  is good and $\Phi^-(D_{i-1},D_i)>0.496$. (To say that
$[D_{i-1},D_{i}]$  is good means that the quantities  $H^-(D_{i-1},D_i)+1$, 
$
\sinh^2(R-D_i)-H^+(D_{i-1},D_i)$, and $\sinh \omega(\tryepsilon/2,D_{i-1})-\sinh
\omega(\tryepsilon/2,D_i)\cdot\sin \theta(\tryepsilon/2,D_i)$ are positive. The smallest
values of these respective quantities, as $i$ ranges from $1$ to $47$,
are $0.75\ldots$, $2.22511\ldots$, and $0.31\ldots$,
and are respectively achieved when $i=1$, $i=47$, and $i=1$.
The
smallest value of $\Phi^-(D_{i-1},D_i)$ is equal to $0.49603\ldots$,
and
is achieved when $i=45$.)
\EndProof

As was mentioned in the introduction,
a group $G$ is
said to be {\it $k$-semifree} for a given positive integer $k$ if each 
subgroup of $G$ having rank at most $k$ is a free product of free
abelian groups.

\Proposition\label{pi-one rank}
Let $M$ be a \willbefinitevolume\ orientable hyperbolic $3$-manifold such that $\pi_1(M)$ is $2$-\willbesemifree. Then
$$\rank\pi_1(M)<1+\willbelambdanought\cdot\vol M,$$
where $\willbelambdanought=156/
\maybeb((\log3)/2)
)=167.781\ldots$.
\EndProposition

\Proof
According to \cite[Corollary 4.2]{acs-surgery}, $\log3$ is a Margulis
number for $M$. 
We will apply Corollary \ref{quote me},
taking $\tryepsilon=\log3$, $R=2\log3+0.15$, and $c=0.496$.
Note that we have
$2\tryepsilon<R<5\tryepsilon/2$, as required for Corollary \ref{quote me}.

Lemma
\ref{compute} asserts that the constant $c$ satisfies
Condition (a) of Corollary \ref{quote me}. The constant $c$ also 
satisfies Condition (b) of that corollary, because $B((\log3)/2)=0.73\ldots>c$.
We have
$(B(R)-\maybeb(\tryepsilon/2))/c=314.62\ldots$, so that Condition (c) of
Corollary \ref{quote me} holds and
$\lfloor(B(R)-\maybeb(\tryepsilon/2))/c\rfloor=314$. The assertion now follows
from Corollary \ref{quote me}.
\EndProof

The transition from Proposition \ref{pi-one rank} to the homology
bounds given by Theorem \ref{main} will involve the following lemma.

\Lemma\label{fuzzwart} 
Let $M$ be a non-compact, finite-volume, orientable hyperbolic
$3$-manifold. Suppose that for some prime $p$ we have
$\dim H_1(M,\FF_p)\ge 3$.
Then
$\vol M>2.848$. 
\EndLemma

\Proof
If $M$ has two cusps, this follows from 
\cite[Theorem 3.6]{agoltwocusp},
which asserts that the smallest volume of an orientable two-cusped
hyperbolic $3$-manifold is $\voct=3.66\ldots$, the volume of a regular
ideal hyperbolic octahedron. If $M$ has more than two
cusps,
it is a standard consequence of Thurston's
hyperbolic Dehn filling theorem \cite[Chapter E]{bp} that there is a hyperbolic manifold $M'$ having
exactly two cusps which can be obtained from $M$ by Dehn filling, and
that $\vol M'<\vol M$. Since $\vol M'\ge\voct$, we have 
$\vol M >\voct$ in this case.

If $M$ has one cusp and $\vol M\le2.848$, then 
\cite[Theorem 1.2]{milley} asserts that $M$ is one of the manifolds
m003, m004, m006, m007, m009, m010, m011, m015, m016 or m017 in the
SnapPea census.
One calculates, using \cite{snappy},
that if $M$ is any of these manifolds then $H_1(M;\FF_p)$ is a direct
sum of $\ZZ$ with a (possibly trivial) cyclic group. Hence
$\dim H_1(M,\FF_p)\le 2$ for every prime $p$.
\EndProof

\Theorem\label{main}
Let $M$ be any \willbefinitevolume\ orientable hyperbolic
$3$-manifold. Then:
\begin{enumerate} 
\item for any prime $p$ we have
$$\dim H_1(M;\FF_p)< \lambda_1\cdot\vol M,$$
where 
$$\lambda_1=\frac1{1.22}+\frac{156}{
b((\log3)/2)}
=168.601\ldots;$$
\item
if $M$ is non-compact, and $p$ is any prime, we have
$$\dim H_1(M;\FF_p)< \lambda_1'\cdot\vol M,$$
where 
$$\lambda_1'=\frac1{2.848}+\frac{156}{
b((\log3)/2)}
=168.132\ldots;$$ and
\item
if $M$ is compact, we have
$$\dim H_1(M;\FF_2)< \lambda_1''\cdot\vol M,$$
where 
$$\lambda_1''=\frac1{3.77}+\frac{156}{
b((\log3)/2)}
=168.046\ldots.$$
\end{enumerate}
\EndTheorem

\Proof
We set $V=\vol M$. Let $p$ be a prime, and set $h_{M,p}=\dim H_1(M;\FF_p)$. We must show:
\Equation\label{first mouse}
h_{M,p}<\lambda_1V;
\EndEquation
\Equation\label{second mouse}
h_{M,p}<\lambda_1'V\text{ if }M\text{ is non-compact; and}
\EndEquation
\Equation\label{third mouse}
h_{M,p}<\lambda_1''V\text{ if }M\text{ is compact and }p=2.
\EndEquation

According to 
\cite[Theorem 1.3]{milley},
we have $\vol M>0.94$ for any \willbefinitevolume\ orientable hyperbolic
$3$-manifold $M$. Hence in the case where 
$h_{M,p}\le10$,
we
have $h_{M,p}<11\cdot V$, which is  stronger than each of the
inequalities (\ref{first mouse})---(\ref{third mouse}).

The rest of the proof will be devoted to the case in which
$h_{M,p}\ge11$.  

Since in particular  $h_{M,p}\ge4$, it
follows from \cite[Lemma 5.2]{acs-surgery} that
$\pi_1(M)$ is  $2$-\willbesemifree.
According to Proposition  
\ref{pi-one rank}, 
we therefore have
$\rank\pi_1(M)<1+\willbelambdanought\cdot V$, 
where $\willbelambdanought=156/
b((\log3)/2)
$.
In particular we have $h_{M,p}<1+\willbelambdanought\cdot V$, which we rewrite as
\Equation\label{get me rewrite}
h_{M,p}<
\bigg(\frac1{V}
+\willbelambdanought\bigg)V.
\EndEquation

Consider the subcase in which
$M$ is compact and $p=2$. Since 
$h_{M,2}\ge11$, it then follows from \cite[Proposition 13.4]{gs} that $V>3.77$, which with
(\ref{get me rewrite}),  gives $h_{M,2}<\lambda_1''V$. This proves
(\ref{third mouse}).

Now consider the subcase in which $M$ is non-compact. Since
$h_{M,p}\ge11>3$, 
Lemma \ref{fuzzwart} gives $V>2.848$. 
With
(\ref{get me rewrite}), this gives $h_{M,p}<\lambda_1'V$, and
(\ref{second mouse}) is proved.

To prove  (\ref{first mouse}), we first notice that in the subcase where $M$ is
non-compact, the asserted inequality follows from (\ref{second mouse}). If $M$ is compact, we use 
\cite[Theorem 1.1]{acs-surgery}, which implies that if $M$ is a 
closed, orientable hyperbolic $3$-manifold with
 $\dim H_1(M;\FF_p)\ge4$ for some prime $p$, then $\vol M>1.22$.
In the notation of the present proof, since $h_{M,p}\ge11>4$, we have
$V>1.22$; with
(\ref{get me rewrite}), this gives $h_{M,p}<\lambda_1V$, and 
(\ref{first mouse}) is proved. 
\EndProof

\bibliographystyle{plain}

\end{document}